\documentclass[paper=a4, USenglish, numbers=noenddot]{scrartcl} 
\usepackage{cmap}		
\usepackage[utf8]{inputenc}	
\usepackage[T1]{fontenc}	
\usepackage[USenglish]{babel}
\usepackage{amsmath}
\usepackage{amssymb}
\usepackage{enumerate}
\usepackage{aliascnt}		
\usepackage{tikz}		
\usetikzlibrary{matrix,arrows,patterns,intersections,calc,decorations.pathmorphing}
\usepackage{calc}
\usepackage{mathtools}		
\usepackage[pdftitle={Rotational component spaces for infinite-type translation surfaces},pdfauthor={Lucien Clavier, Anja Randecker, Chenxi Wu},pdfborder={0 0 0}]{hyperref}   
\usepackage[figure]{hypcap}	

\usepackage{amsthm}

\newtheoremstyle{break}
 {} 
 {} 
 {\itshape} 
 {} 
 {\bfseries} 
 {} 
 {\newline} 
 {\thmname{#1}\thmnumber{ #2}\thmnote{ (#3)}} 
 
\newtheoremstyle{breakdef}
 {} 
 {} 
 {} 
 {} 
 {\bfseries} 
 {} 
 {\newline} 
 {\thmname{#1}\thmnumber{ #2}\thmnote{ (#3)}} 
 
\newtheoremstyle{remark}
 {} 
 {} 
 {} 
 {} 
 {\itshape} 
 {.} 
 {0.5em} 
 {\thmname{#1}\thmnumber{ #2}\thmnote{ {\normalfont (}#3{\normalfont )}}} 

\theoremstyle{breakdef}
\newtheorem{definition}{Definition}[section]

\theoremstyle{remark}

\newaliascnt{rem}{definition}  

\aliascntresetthe{rem}

\newaliascnt{exa}{definition}  
\newtheorem{exa}[exa]{Example}
\aliascntresetthe{exa} 

\newaliascnt{cor}{definition}  
\newtheorem{cor}[cor]{Corollary}
\aliascntresetthe{cor}  

\newaliascnt{lem}{definition}  
\newtheorem{lem}[lem]{Lemma}
\aliascntresetthe{lem}

\theoremstyle{break}

\newaliascnt{prop}{definition}  
\newtheorem{prop}[prop]{Proposition}
\aliascntresetthe{prop}

\newaliascnt{thm}{definition}  
\newtheorem{thm}[thm]{Theorem}
\aliascntresetthe{thm}

\numberwithin{figure}{section}

\newcommand{\LC}{\mathcal{L}} 
\newcommand{\LT}{\widetilde{\mathcal{L}}}

\newcommand{\RR}{\mathbb{R}} 
\newcommand{\ZZ}{\mathbb{Z}}

\author{Lucien Clavier, Anja Randecker, Chenxi Wu}
\title{Rotational component spaces for infinite-type translation surfaces}
\date{\today}

\begin{document}
\renewcommand{\sectionautorefname}{Section}
 
\maketitle

\begin{abstract}
 Finite translation surfaces can be classified by the order of their singularities. When generalizing to infinite translation surfaces, however, the notion of order of a singularity is no longer well-defined and has to be replaced by new concepts.
 
 This article discusses the nature of two such concepts, recently introduced by Bowman and Valdez: linear approaches and rotational components. We show that there is a large flexibility in the spaces of rotational components and even more in the spaces of linear approaches. In particular, we prove that every finite topological space arises as space of rotational components. However, this space will still not contain enough information to describe an infinite translation surface. We showcase this through an uncountable family with the same space of rotational components but different spaces of linear approaches.

 Additionally, we study several known and new examples to illustrate the concept of linear approaches and rotational components.
\end{abstract}

Perhaps the first appearance of translation surfaces was in the context of dynamics in~\cite{fox_kershner_36}, more specifically in the context of billiards. The unfolding construction for billiards leads to the most visual way to describe translation surfaces: they can be obtained by gluing polygons in the Euclidean plane along parallel edges of the same length.

Since the seminal works of Earle, Gardiner, Kerckhoff, Masur, Smillie, Veech, Vorobets, and many more, starting in the 80's, the interest in this field of research has grown and different aspects have been studied. One aspect is the study of the geometry of moduli spaces of finite translation surfaces which can be broken down into strata (see \cite{zorich2006flat} for a survey). For these strata, for example, the number and type of connected components is known since \cite{kontsevich2003connected}. Nowadays, dynamics on strata is a rich theory and its examination has recently lead to the breakthrough result on orbit closures in \cite{eskin2013isolation}.

Within the last few years, translation surfaces of infinite topological type have come into the focus of research. Many results have been obtained for infinite coverings of finite translation surfaces (e.g.\ \cite{delecroix_hubert_lelievre_11}, \cite{hooper_hubert_weiss_13}, and \cite{fraczek_ulcigrai_14}). Furthermore, unfamiliar behaviours of the dynamics in the infinite case have been found in \cite{hooper_10}, \cite{trevino_12}, and \cite{lindsey_trevino_14}, for~example.

One of the prominent aims in this young area of research is to extend the classical theory of dynamics on strata to the world of infinite-type surfaces.
In the finite-type case, neighborhoods of singularities are simply covers of the punctured disk, and we can use the degree of these covers to define strata. 
In order to define an analog to strata in the general case, it is natural to first study
which properties neighborhoods of general singularities can have.
This was first undertaken in \cite{bowman_valdez2013wild}, where the authors studied the germs of geodesic rays starting at singularities, called \emph{linear approaches} (see \autoref{def:linear_approaches}). 
All linear approaches of a singularity $\sigma$ together form a space $\LC(\sigma)$.
This space can be seen as the tangent space at the singularity $\sigma$ and so it helps to extend the unit tangent bundle to the singularities of the translation surface.

If two linear approaches differ by a rotation centered at the singularity, we say they are in the same \emph{rotational component} (see \autoref{def:rotational_component}).
Being in the same rotational component constitutes an equivalence relation~$\sim$ on $\LC(\sigma)$.
For regular points and for the well-known kinds of singularities, there is only one rotational component.

The rotational components are for example useful as a tool to study the dynamics of a translation surface. This is done in a book in preparation by Delecroix, Hubert, and Valdez and in an unpublished article by Rafi and Randecker.

One can also study the set $\LT(\sigma) = \LC(\sigma)/\sim$ of all rotational components, equipped with the quotient topology. This space captures how rotational components “fit together”.
When looking at known examples of singularities of translation surfaces and their spaces of rotational components, it seemed reasonable to conjecture that they have properties such as being locally compact or having Lebesgue covering dimension $1$. It could have been even possible that a space of rotational components is always Hausdorff, inheriting this property from the corresponding space of linear approaches.

Quite in contrast, we show in the present paper that there is a large variety of spaces that can occur as spaces of rotational components. In particular, a space of rotational components does not need to be Hausdorff and the Lebesgue covering dimension can be arbitrarily high. This follows from \autoref{thm:finite_top}.

{
\renewcommand{\thethm}{\ref{thm:finite_top}}
\begin{thm}[Every finite space occurs as $\LT(\sigma)$]
 Let $Y$ be a non-empty topological space of finite cardinality. Then there exists a compact translation surface $(X, \mathcal{A})$ with a wild singularity $\sigma$ so that the space $\LT(\sigma)$ of rotational components carries the same topology as $Y$.
\end{thm}
}

This theorem implies, moreover, that $\LT(\sigma)$ can be non-$T_0$, $T_0$ but not $T_1$, $T_1$ but not $T_2$, or $T_2$.
Examples in \autoref{sec:examples} illustrate these different possibilities and \autoref{exa:non-locally_compact} gives an instance of a non-locally compact space of rotational components. On the other hand, \autoref{exa:star_decoration} shows that the space of rotational components can also be homeomorphic to $S^1$.

In general, it is unlikely that a space of rotational components inherits properties of the corresponding space of linear approaches. This is because we lose a large amount of information when passing to the quotient as we show in \autoref{thm:info_loss}.

{
\renewcommand{\thethm}{\ref{thm:info_loss}}
\begin{thm}[$\LT(\sigma)$ does not determine $\LC(\sigma)$]
 There are uncountably many translation surfaces $(X_r, \mathcal{A}_r)$ with exactly one singularity~$\sigma_r$, with pairwise non-homeomorphic spaces of linear approaches $\LC(\sigma_r)$, but with homeomorphic spaces of rotational components $\LT(\sigma_r)$.
\end{thm}
}

Furthermore, we investigate relations between the behaviour of the rotational components as points in $\LT(\sigma)$ and a natural metric that we can put on each rotational component. In \autoref{prop:biinfinite_means_open} and \autoref{prop:open_means_infinite_length}, we prove that the length of a rotational component is related to whether it defines an open point in $\LT(\sigma)$. With similar arguments as in these proofs, we can also show a constraint on the spaces of rotational components: In \autoref{cor:separable}, we verify that $\LT(\sigma)$ has to be separable.

\bigskip
\noindent\textbf{Acknowledgements.}
The idea for this project was developed during the semester program “Low-dimensional Topology, Geometry, and Dynamics” at the Institute for Computational and Experimental Research in Mathematics in Providence.
Most of this work was done while the first and the third author were visiting the Karlsruhe Institute of Technology, funded by a Visiting Researcher Scholarship of the Karlsruhe House of Young Scientists.
The second author was partially supported within the project “Dynamik unendlicher Translationsflächen” in the “Juniorprofessoren-Programm” of the Ministry of Science, Research and the Arts of Baden-Wuerttemberg.

\section{Definitions}\label{sec:def}

A \emph{translation surface} $(X,\mathcal{A})$ is a two-dimensional connected manifold $X$ equipped with a translation structure $\mathcal{A}$, i.e.\ a maximal atlas whose transition functions are translations.
We can pull back the Euclidean metric to the translation surface via the charts.
When passing to the metric completion we possibly have additional points which are called \emph{singularities}.

For a long time, only \emph{finite translation surfaces} were studied. A finite translation surface comes from a compact Riemann surface with a holomorphic \mbox{$1$--form}. The \mbox{$1$--form} induces a translation structure on the surface without the zeros of the \mbox{$1$--form}.
Here, the singularities are points where curvature is concentrated and which are cone points with cone angle $2\pi k$ for a $k\geq 2$.
These points are called \emph{cone angle singularities} of multiplicity~$k$. They have a punctured neighborhood which is a cyclic translation covering of a once-punctured disk of degree $k$.

We say that a translation surface is of \emph{infinite type} if its metric completion has singularities that are not cone angle singularities.
By small abuse of notation we often identify a translation surface with its metric completion.
For example we say a translation surface $(X, \mathcal{A})$ is \emph{compact} if the metric completion $\overline{X}$ is compact.

An infinite-type translation surface can have singularities for which the corresponding cyclic translation covering is infinite (these are called \emph{infinite angle singularities}). Even more, an infinite-type translation surface can also have singularities that do not have a neighborhood which is a cyclic translation covering of a once-punctured disk.
The latter singularities are called \emph{wild singularities}. They incorporate most of the complication inherent to infinite-type translation surfaces.
Therefore, even when statements are true for all types of singularities, we will mainly be interested in the case of wild singularities.

In \cite{bowman_valdez2013wild}, Joshua Bowman and Ferrán Valdez initiated the classification of wild singularities by studying their spaces of linear approaches.
We recall the definitions that we will use here.

\begin{definition}[Space of linear approaches]\label{def:linear_approaches}
 Let $(X,\mathcal{A})$ be a translation surface. For every $\epsilon > 0$, we define the space
 \begin{equation*}
  \LC^\epsilon (X) \coloneqq \left\{ \gamma \colon (0,\epsilon) \to X : \gamma \text{ is a geodesic curve} \right\}
  .
 \end{equation*}
 On $\bigsqcup\limits_{\epsilon > 0} \LC^\epsilon (X)$, we have the following equivalence relation $R$: $\gamma_1 \in \LC^\epsilon (X)$ and $\gamma_2 \in \LC^{\epsilon'} (X)$ are called $R$--equivalent if $\gamma_1 (t) = \gamma_2 (t)$ for all $t\in (0,\min\{\epsilon, \epsilon' \})$.
 The space
 \begin{equation*}
  \LC (X) \coloneqq \bigsqcup\limits_{\epsilon > 0} \LC^\epsilon (X) / R
 \end{equation*}
 is called \emph{space of linear approaches of $X$} and the $R$--equivalence class $[\gamma]$ of $\gamma\in \LC^{\epsilon}(X)$ is called \emph{linear approach}.
\end{definition}

For a point $x\in \overline{X}$, we define the subspace $\LC^\epsilon(x) \coloneqq \left\{ \gamma \in \LC^\epsilon(X) : \lim\limits_{t \to 0} \gamma(t) = x \right\}$ and $\LC(x)$ as before. The $R$--equivalence class $[\gamma]$ of $\gamma\in \LC^{\epsilon}(x)$ is called \emph{linear approach to the point~$x$}.

Before specifying the topology on $\LC(X)$ and $\LC(x)$, we define rotational components as classes of linear approaches. We also define the translation structure on a rotational component.

Following \cite{bowman_valdez2013wild}, by an \emph{angular sector} we mean a triple $(I,c,i_c)$ such that $I$ is a non-empty generalized interval (i.e.\ a connected subset of $\mathbb{R}$), $c\in \mathbb{R}$ and $i_c$ is an isometry from $\{x+\mathrm{i}y : x<c,y\in I\} \subset \mathbb{C}$ with metric defined by $e^zdz$ to $X$.
Given any angular sector $(I,c,i_c)$ with $z \coloneqq \lim_{x \to -\infty} \{ i_c(x+\mathrm{i}y ) \}$, there is a map $f_{(I,c,i_c)} \colon I\rightarrow \LC(z)$ which sends an element $y\in I$ to the linear approach that is represented by the ray $i_c(t+\mathrm{i}y)$ for~$t<c$.

\begin{definition}[Rotational component]\label{def:rotational_component}
 We define an equivalence relation $\sim$ on $\LC(X)$ as follows: two linear approaches $[\gamma_1]$ and $[\gamma_2]$ are $\sim$--equivalent if there is an angular sector $(I,c,i_c)$ and $i_1, i_2 \in I$ such that $f_{(I,c,i_c)}(i_1)=[\gamma_1]$ and $f_{(I,c,i_c)}(i_2)=[\gamma_2]$.

 An equivalence class under $\sim$ is called a \emph{rotational component}. We write $\overline{[\gamma]}$ for the equivalence class of a linear approach $[\gamma]$.
\end{definition}

\begin{definition}[Translation structure on rotational components]\label{def:translation_struct}
 The image of $f_{(I,c,i_c)}$ is always contained in a single rotational component, hence the map $f_{(I,c,i_c)}$ gives a chart of a translation structure with boundary for every rotational component that contains more than one linear approach.
\end{definition}

We now describe the topology with which we endow $\LC(X)$.

\begin{definition}[Topology on $\LC(X)$] \label{def:topology_on_LCX}
 We equip $\LC^\epsilon(X)$ with the topology induced by the metric
 \begin{equation*}
  d_\epsilon(\gamma_1,\gamma_2)=\sup_{t\in(0,\epsilon)}d(\gamma_1(t), \\ \gamma_2(t)) .
 \end{equation*}

 For every $\epsilon > 0$, we can embed $\LC^{\epsilon}(X)$ in $\LC(X)$ so that we obtain a direct system. Then we define the topology on $\LC(X)$ as the colimit topology with respect to the embeddings $\LC^\epsilon(X) \hookrightarrow \LC(X)$.
 On $\LC^\epsilon(x)$ and $\LC(x)$ for $x\in \overline{X}$, we 
 consider the subspace topology as subspaces of $\LC^\epsilon(X)$ and $\LC(X)$.
\end{definition}

This definition of the topology on $\LC(X)$ is quite abstract so we recall a characterization from \cite[Proposition 2.3]{bowman_valdez2013wild} which is easier to use in our concrete examples.
For every $x\in X$, $r>0$, and $t>0$, the~set
\begin{equation*}
 \widetilde{B} (x,r)^t=\{[\gamma]: \gamma\in\LC^{\epsilon}(X) \text{ for some } \epsilon>t,d(\gamma(t),x)<r\}
\end{equation*}
is an open set.
In fact, the collection $\{\widetilde{B}(x,r)^t\}$ forms a subbasis for the previously described topology on $\LC(X)$.

In particular, this leads to a characterization of the topology on the space $\LC(x)$ of linear approaches of $x \in X$ which we will often use in proofs and computations.

\begin{lem}[Subbasis for the topology on $\LC(x)$]\label{lem:subbasis_LCpoint}
 Let $(X, \mathcal{A})$ be a translation surface.
 For every geodesic curve $\gamma \in \LC^\epsilon (x)$, $r > 0$, and $t \in (0, \epsilon)$, the set
 \begin{equation*}
  B(\gamma, t, r) \coloneqq \widetilde{B}( \gamma(t), r) ^t \cap \LC(x)
  = \left\{ [\gamma']: \gamma'\in\LC^{\epsilon'}(x) \text{ for some } \epsilon'>t, \ d(\gamma(t), \gamma'(t) )<r \right\}
 \end{equation*}
 is an open neighborhood of $[\gamma]$ in $\LC(x)$.
 In fact, the collection $\{B (\gamma, t, r) \}$ forms a subbasis for the topology on $\LC(x)$ as a subspace of $\LC(X)$.

\begin{proof}
 The first statement is clear by the definition of the topology on $\LC(x)$ and by $[\gamma] \in B(\gamma, t, r)$.

 To prove the second statement, let $x'\in X$, $r'>0$, and $t'>0$. For every linear approach $[\gamma] \in \widetilde{B}(x',r')^{t'} \cap \LC(x)$, we specify a neighborhood from the given collection which is contained in $\widetilde{B}(x',r')^{t'} \cap \LC(x)$:
 Let $\gamma$ be a representative of $[\gamma]$ and define $t\coloneqq t'$ and $r > 0$ small enough so that $B(\gamma(t), r) \subseteq B(x', r')$. Then we have
 \begin{equation*}
  B( \gamma, t, r) \subseteq \widetilde{B}(x', r')^{t'} \cap \LC(x)
 \end{equation*}
 which shows that $\{B (\gamma, t, r) \}$ is a subbasis for the topology on $\LC(x)$.
\end{proof}
\end{lem}

As a direct consequence of the definitions, we have the following very useful lemma (see \cite[Corollary 2.2]{bowman_valdez2013wild}).

\begin{lem}[Directions of linear approaches are varying continuously]\label{lem:dir_continuous}
 Let $(X, \mathcal{A})$ be a translation surface. Then the map $\LC(X)\to S^1$ that associates to a linear approach its direction is continuous.
\end{lem}

The space of rotational components can now also be provided with a topology. 

\begin{definition}[Space of rotational components]
 We define the \emph{space of rotational components} as $\LT(X) \coloneqq \LC(X)/\sim$ endowed with the quotient topology.

 Furthermore, we define $\LT(x)$ as $\LT(x)\coloneqq\LC(x)/\sim$, and give it the topology induced by the inclusion into $\LT(X)$.
\end{definition}

\section{Examples}\label{sec:examples}

Many interesting examples of spaces of rotational components arise from starting with a translation surface with desirable properties, cutting it along geodesic segments and regluing the segments in a different way.
By performing this operation carefully, we are able to control some of the properties of the resulting translation surface, for example compactness or finite area.
This method is also used in \cite{przytycki_schmithuesen_valdez_11} to construct a translation surface with a given Veech group.

Let $(X_1, \mathcal{A}_1)$, $(X_2, \mathcal{A}_2)$ be translation surfaces (not necessarily different).
A \emph{slit} $m$ on $(X_1, \mathcal{A}_1)$ is an open geodesic segment in $X_1$ with an orientation.
When we consider the metric completion of $X_1 \setminus m$ equipped with the path-length metric, we obtain two copies $m'$ and $m''$ of $m$, with endpoints identified.
If $m_1$ is a slit on $(X_1, \mathcal{A}_1)$ and $m_2$ is a slit on $(X_2, \mathcal{A}_2)$ with the same holonomy vectors then we can “glue the slits together”.
For instance, we can identify $m_1'$ with $m_2''$ or $m_1''$ with $m_2'$ by a translation and we will describe this visually as “gluing the upper part of $m_1$ to the lower part of $m_2$”.
If we glue both $m_1'$ with $m_2''$ and $m_1''$ with $m_2'$, we just say that we glue $m_1$ and $m_2$ together.
Another possibility is to consider one slit and divide it into several segments which we reglue in an order that differs from the original order.

As a first example, we decorate the Euclidean plane using this construction.
The space of rotational components of the resulting surface is just a singleton, so the topology is already determined.

\begin{exa}[Harmonic series decoration] \label{exa:harmonic_series_decoration}
 Consider the Euclidean plane $\mathbb{R}^2$ with a horizontal slit starting from $(0,0)$ going to the right forever.
 We choose a sequence $(a_n)$ in $\mathbb{R}_+$ such that $\lim_{n\rightarrow\infty}a_n=0$ and $\sum_{n=1}^\infty a_n=\infty$, e.g.\ the harmonic series.
 We divide the upper and the lower part of the slit into segments. For the upper part the segments are $I_1$ of length $a_1$, $I_2$ of length $a_2$, and so on inductively, for the lower part the segments are $J_2$ of length $a_2$, $J_1$ of length $a_1$, then $J_4$, $J_3$, $J_6$, and so on (see \autoref{fig:harmonic_series_decoration}).
 Then we glue each $I_n$ with $J_n$.

 \begin{figure}[bthp]
 \begin{center}
 \begin{tikzpicture}[xscale=4, yscale=2]
  \draw (0,0) .. node[above] {$I_1$} controls +(20:0.2cm) and +(-180:0.3cm) .. (1,0.1)
   -- node[above] {$I_2$} (1.5,0.1)
   -- node[above] {$I_3$} (1.833,0.1)
   -- node[above] {$I_4$} (2.083,0.1)
   -- node[above] {$I_5$} (2.283,0.1)
   -- node[above] {$I_6$} (2.45,0.1) -- (2.5,0.1);
  \draw[densely dotted] (2.5,0.1) -- (2.6,0.1);
  \draw (1,0.12) -- (1,0.08);
  \draw (1.5,0.12) -- (1.5,0.08);
  \draw (1.833,0.12) -- (1.833,0.08);
  \draw (2.083,0.12) -- (2.083,0.08);
  \draw (2.283,0.12) -- (2.283,0.08);
  \draw (2.45,0.12) -- (2.45,0.08);

  \draw (0,0) .. node[below] {$J_2$} controls +(-20:0.2cm) and +(175:0.1cm) .. (0.5,-0.08)
   .. node[below] {$J_1$} controls +(-5:0.2cm) and +(180:0.5cm) .. (1.5,-0.1)
   -- node[below] {$J_4$} (1.75,-0.1)
   -- node[below] {$J_3$} (2.083,-0.1)
   -- node[below] {$J_6$} (2.25,-0.1)
   -- node[below] {$J_5$} (2.45,-0.1) -- (2.5,-0.1);
  \draw[densely dotted] (2.5,-0.1) -- (2.6,-0.1);
  \draw (0.5,-0.06) -- (0.5,-0.1);
  \draw (1.5,-0.12) -- (1.5,-0.08);
  \draw (1.75,-0.12) -- (1.75,-0.08);
  \draw (2.083,-0.12) -- (2.083,-0.08);
  \draw (2.25,-0.12) -- (2.25,-0.08);
  \draw (2.45,-0.12) -- (2.45,-0.08);
 \end{tikzpicture}
 \label{fig:harmonic_series_decoration}
 \caption{Harmonic series decoration: segments $I_i$ and $J_i$ are glued.}
 \end{center}
 \end{figure}
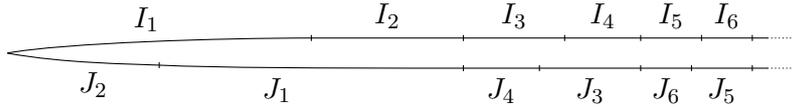

 All the endpoints of the segments are identified via the gluings, so there is only one singularity $\sigma$, which is obviously not a cone angle singularity.
 Also, $\sigma$ cannot be an infinite angle singularity as for every $\epsilon > 0$ there exists a saddle connection of length smaller than~$\epsilon$ starting and ending at $\sigma$. 
 Thus, $\sigma$ is a wild singularity. It has exactly one rotational component which is isometric to $\RR$ as a class of linear approaches (see \autoref{def:translation_struct}).
 So $\LT(\sigma)$ is the one-point space equipped with the unique possible topology.
\end{exa}

The following is a refinement of the previous example with non-trivial space of rotational components.

\begin{exa}[Geometric series decoration]\label{exa:geometric_series_decoration}
 Consider a translation surface with a slit of finite length. We choose a sequence $(a_n)$ in $\mathbb{R}_+$ for which the corresponding series is converging to the length of the slit, e.g.\ a geometric series.
 We divide the upper and the lower part of the slit into segments of length $a_n$ and identify the segments of same index.
 We will describe two different procedures how to do the construction.

 \begin{enumerate}
  \item First we can choose to assign a segment $I_1$ (resp.\ $J_1$) of length $a_1$ at the left (resp.\ right) of the upper (resp.\ lower) part.
  Then assign a segment $I_2$ (resp.\ $J_2$) of length~$a_2$ at the right of $I_1$ (resp.\ at the left of $J_1$), and so on inductively (see \autoref{fig:geometric_series_decoration_standard}).
  
  In this case we obtain one singularity $\sigma$ with two rotational components, both isometric to $(0,\infty)$. 
  Let $[\gamma_1]$ be the vertical linear approach starting on the far right of the slit and going upward. Its corresponding rotational component shall be~$\overline{[\gamma_1]}$.
  Then every vertical linear approach starting at the right endpoint of $I_{2n}$ for $n\geq 1$ and going upward is contained in the other rotational component, called $\overline{[\gamma_2]}$. 
  Thus there is a sequence of vertical linear approaches belonging to $\overline{[\gamma_2]}$ converging to~$[\gamma_1]$ in the topology on $\LC(\sigma)$.
  We can conclude that any open set in $\LT(\sigma)$ containing~$\overline{[\gamma_1]}$ also contains $\overline{[\gamma_2]}$, and by symmetry we see that the topology of $\LT(\sigma)$ is
  \begin{equation*}
   \left\{ \emptyset , \{\overline{[\gamma_1]}, \overline{[\gamma_2]}\} \right\}.
  \end{equation*}

  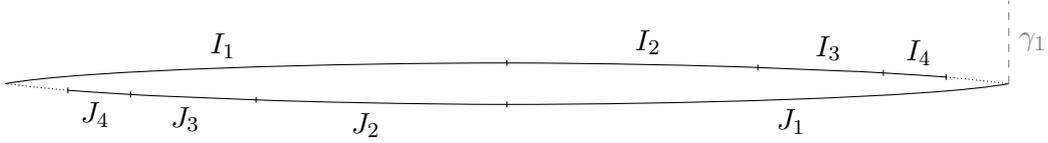
\begin{figure}[bhtp]
  \begin{center}
  \begin{tikzpicture}[xscale=1.1, yscale=0.55]
   \draw (0,0) .. node[above] {$I_1$} controls +(20:1cm) and +(-180:2cm) .. (6,0.5);
   \draw (6,0.57) -- (6,0.43);
   \draw (6,0.5) .. node[above] {$I_2$} controls +(0:1cm) and +(175:0.5cm) .. (9,0.39);
   \draw (9,0.46) -- (9,0.32);
   \draw (9,0.39) .. node[above] {$I_3$} controls +(-5:0.5cm) and +(172:0.25cm) .. (10.5,0.26);
   \draw (10.5,0.33) -- (10.5,0.19);
   \draw (10.5,0.26) .. node[above] {$I_4$} controls +(-8:0.25cm) and +(169:0.12cm) .. (11.25,0.16);
   \draw (11.25,0.23) -- (11.25,0.09);
   \draw[densely dotted] (11.25,0.16) .. controls +(-11:0.12cm) and +(160:0.06cm) .. (12,0);
   \draw (12,0) .. node[below] {$J_1$} controls +(-160:1cm) and +(0:2cm) .. (6,-0.5);
   \draw (6,-0.57) -- (6,-0.43);
   \draw (6,-0.5) .. node[below] {$J_2$} controls +(180:1cm) and +(-5:0.5cm) .. (3,-0.39);
   \draw (3,-0.46) -- (3,-0.32);
   \draw (3,-0.39) .. node[below] {$J_3$} controls +(175:0.5cm) and +(-8:0.25cm) .. (1.5,-0.26);
   \draw (1.5,-0.33) -- (1.5,-0.19);
   \draw (1.5,-0.26) .. node[below] {$J_4$} controls +(172:0.25cm) and +(-11:0.12cm) .. (0.75,-0.16);
   \draw (0.75,-0.23) -- (0.75,-0.09);
   \draw[densely dotted] (0.75,-0.16) .. controls +(169:0.12cm) and +(-20:0.06cm) .. (0,0);
   
   \draw[color=gray,dashed] (12,0) -- node[right]{$\gamma_1$} (12,2);
  \end{tikzpicture}
  \label{fig:geometric_series_decoration_standard}
  \caption{Geometric series decoration with segments going from left to right on the upper part and going from right to left on the lower part.}
  \end{center}
  \end{figure}

  \item A different topology on a two-element space of rotational components can be obtained in the following way.
  Similar to \autoref{exa:harmonic_series_decoration}, we first assign a segment $I_1$ of length $a_1$ at the left of the upper part of the slit and assign successive segments~$I_2$ of length $a_2$ and so on inductively on its right.
  For the lower part, we assign a segment $J_2$ of length $a_2$ at the left, a segment $J_1$ of length $a_1$ on its right, then $J_4$ of length $a_4$, $J_3$ of length $a_3$, and so on.
  Then we glue the segments of the same length (see \autoref{fig:geometric_series_decoration_modified}).

  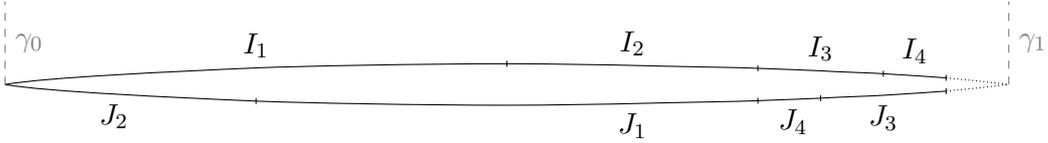
\begin{figure}
  \begin{center}
  \begin{tikzpicture}[xscale=1.1, yscale=0.55]
   \draw (0,0) .. controls +(20:0.5cm) and +(-175:1cm) .. (3,0.4) node[above] {$I_1$} .. controls +(5:0.5cm) and +(-180:1cm) .. (6,0.5);
   \draw (6,0.57) -- (6,0.43);
   \draw (6,0.5) .. controls +(0:0.5cm) and +(177:0.25cm) .. (7.5,0.46) node[above] {$I_2$} .. controls +(-3:0.5cm) and +(175:0.25cm) .. (9,0.39);
   \draw (9,0.46) -- (9,0.32);
   \draw (9,0.39) .. controls +(-5:0.25cm) and +(173:0.125cm) .. (9.75,0.33) node[above] {$I_3$} .. controls +(-7:0.25cm) and +(172:0.125cm) .. (10.5,0.26);
   \draw (10.5,0.33) -- (10.5,0.19);
   \draw (10.5,0.26) .. controls +(-8:0.25cm) and +(169:0.12cm) .. (11.25,0.16);
   \draw (10.875,0.27) node[above] {$I_4$};
   \draw (11.25,0.23) -- (11.25,0.09);
   \draw[densely dotted] (11.25,0.16) .. controls +(-11:0.12cm) and +(160:0.06cm) .. (12,0);
   \draw (0,0) .. node[below] {$J_2$} controls +(-20:0.5cm) and +(175:1cm) .. (3,-0.4);
   \draw (3,-0.47) -- (3,-0.33);
   \draw (3,-0.4) .. controls +(-5:0.5cm) and +(180:1cm) .. (6,-0.5) .. controls +(0:0.5cm) and +(-177:0.25cm) .. (7.5,-0.46) node[below] {$J_1$} .. controls +(3:0.5cm) and +(-175:0.25cm) .. (9,-0.39);
   \draw (9,-0.46) -- (9,-0.32);
   \draw (9,-0.39) .. node[below] {$J_4$} controls +(5:0.25cm) and +(-173:0.125cm) .. (9.75,-0.33);
   \draw (9.75,-0.4) -- (9.75,-0.26);
   \draw (9.75,-0.33) .. controls +(7:0.25cm) and +(-172:0.125cm) .. (10.5,-0.26) node[below] {$J_3$} .. controls +(8:0.25cm) and +(-169:0.12cm) .. (11.25,-0.16);
   \draw (11.25,-0.23) -- (11.25,-0.09);
   \draw[densely dotted] (11.25,-0.16) .. controls +(11:0.12cm) and +(-160:0.06cm) .. (12,0);
   
   \draw[color=gray,dashed] (0,0) -- node[right]{$\gamma_0$} (0,2);
   \draw[color=gray,dashed] (12,0) -- node[right]{$\gamma_1$} (12,2);
  \end{tikzpicture}
  \label{fig:geometric_series_decoration_modified}
  \caption{Geometric series decoration with segments going in the same direction on the upper and the lower part of the slit.}
  \end{center}
  \end{figure}

  Again we have one singularity $\sigma$ with two rotational components.
  For the vertical linear approach $[\gamma_0]$ that is starting on the left endpoint of the slit and going upward, the rotational component $\overline{[\gamma_0]}$ is isometric to $(-\infty, \infty)$.
  On the other hand, for the vertical linear approach $[\gamma_1]$ that is starting on the right endpoint of the slit and going upward, the rotational component $\overline{[\gamma_1]}$ has length $2\pi$.
  With a similar argument as in (i) we can see that there are linear approaches in $\overline{[\gamma_0]}$ converging to a linear approach in $\overline{[\gamma_1]}$. However, there is no sequence in $\overline{[\gamma_1]}$ converging to a linear approach in $\overline{[\gamma_0]}$. 
  Therefore, the topology of $\LT(\sigma)$ is
  \begin{equation*}
   \left\{ \emptyset, \{\overline{[\gamma_0]}\}, \{\overline{[\gamma_0]}, \overline{[\gamma_1]} \} \right\}.
  \end{equation*}
 \end{enumerate}
\end{exa}

The following example can be found in \cite{bowman2012finiteness} and will be very useful in the proof of \autoref{thm:info_loss}.

\begin{exa}[Stack of boxes] \label{exa:stack_of_boxes}
 Choose two sequences $H = \{h_n\}$ and $W = \{w_n\}$ of positive real numbers so that $H$ is bounded away from $0$ (for instance $h_n\geq 1$), $W$ is strictly decreasing, and $\sum_{n=1}^\infty h_nw_n<\infty$.
 Now, an infinite sequence of rectangles with height $h_n$ and width $w_n$ is stacked on top of each other as shown in \autoref{fig:stack_of_boxes}. The bottom of the first rectangle is divided into segments of length $w_n-w_{n+1}$ for every $n\geq 1$. By gluing the opposite sides, all corners of the rectangles and the points $A_i$ are identified. Hence, this surface has only one singularity with exactly one rotational component which is isometric to~$(0,\infty)$.

 We can specify this isometry by starting with a linear approach which starts at the left lower corner $A_0$ of the first rectangle and stays in the first rectangle for some time.
 When rotating this linear approach counterclockwise around $A_0$, we obtain linear approaches that start at $A_1$, then at $B_1, C_1, D_1, A_2, B_2, \ldots$ When we do the rotation clockwise we cannot rotate further than to an almost horizontal linear approach which corresponds to a number almost $0$.

 \begin{figure}
 \begin{center}
 \begin{tikzpicture}
  \draw[densely dotted] (0,0) -- (1,0);
  \draw (1,0) -- (1.7,0) -- (3,0) -- (4.3,0) -- (5,0) -- (5,1) -- (0,1) -- (0,0);
  \draw (4.3,1) -- (4.3,2.2) -- (0,2.2) -- (0,1);
  \draw (3,2.2) -- (3,4.2) -- (0,4.2) -- (0,2.2);
  \draw (1.7,4.2) -- (1.7,5) -- (0,5) -- (0,4.2);
  \draw (1,5) -- (1,7.2) -- (0,7.2) -- (0,5);
  \draw[densely dotted] (0.8,7.2) -- (0.8,7.4);
  \draw[densely dotted] (0,7.2) -- (0,7.4);
  \draw[fill] (0,0) circle (1pt) node[below] {$A_0$};
  \draw[fill] (1,0) circle (1pt) node[below] {$A_5$};
  \draw[fill] (1.7,0) circle (1pt) node[below] {$A_4$};
  \draw[fill] (3,0) circle (1pt) node[below] {$A_3$};
  \draw[fill] (4.3,0) circle (1pt) node[below] {$A_2$};
  \draw[fill] (5,0) circle (1pt) node[below] {$A_1$};
  \draw[fill] (5,1) circle (1pt) node[right] {$B_1$};
  \draw[fill] (4.3,2.2) circle (1pt) node[right] {$B_2$};
  \draw[fill] (3,4.2) circle (1pt) node[right] {$B_3$};
  \draw[fill] (1.7,5) circle (1pt) node[right] {$B_4$};
  \draw[fill] (1,7.2) circle (1pt) node[right] {$B_5$};
  \draw[fill] (0,1) circle (1pt) node[left] {$C_1$};
  \draw[fill] (0,2.2) circle (1pt) node[left] {$C_2$};
  \draw[fill] (0,4.2) circle (1pt) node[left] {$C_3$};
  \draw[fill] (0,5) circle (1pt) node[left] {$C_4$};
  \draw[fill] (0,7.2) circle (1pt) node[left] {$C_5$};
  \draw[fill] (4.3,1) circle (1pt) node[above right] {$D_1$};
  \draw[fill] (3,2.2) circle (1pt) node[above right] {$D_2$};
  \draw[fill] (1.7,4.2) circle (1pt) node[above right] {$D_3$};
  \draw[fill] (1,5) circle (1pt) node[above right] {$D_4$};
 \end{tikzpicture}
 \label{fig:stack_of_boxes}
 \caption{Stack of boxes: opposite sides are glued.}
 \end{center}
 \end{figure}
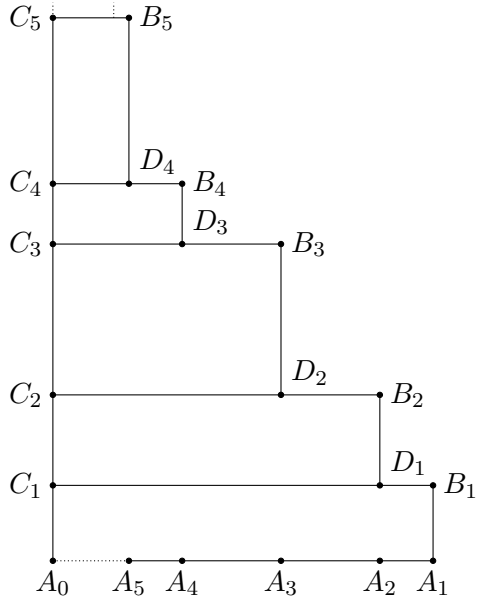
\end{exa}

It is worth noting that the space $\LC(\sigma)$ is Hausdorff in general (see \cite[Corollary 2.9]{bowman_valdez2013wild}). On the other hand, we have the following proposition.

\begin{prop}[$\LC(\sigma)$ is in general not metrizable]
 For the singularity $\sigma$ of the stack of boxes, the space $\LC(\sigma)$ of linear approaches is not T\textsubscript{3}, hence not metrizable.

\begin{proof}
 For a fixed $\epsilon' > 0$, we define the set
 \begin{equation*}
  F=\LC(\sigma) \setminus \Bigl(\bigcup_{\epsilon>\epsilon'} \LC^{\epsilon}(\sigma)\Bigr) .
 \end{equation*}
 Here the spaces $\LC^{\epsilon}(\sigma)$ are identified with their images under the embedding \mbox{$\LC^{\epsilon}(\sigma)\rightarrow\LC(\sigma)$} defined in \autoref{def:topology_on_LCX}.
 The set $F$ is closed as \cite[Corollary 2.5]{bowman_valdez2013wild} shows that its complement is open in $\LC(\sigma)$.

 For $\epsilon'$ small enough, the linear approach $[\gamma]$ starting in the lower right corner $A_1$ and going left toward $A_2$ is not contained in $F$.
 On the other hand, given a point $\gamma(t)$ and an $r \in (0,t)$, there exists an $n\in \mathbb{N}$ such that the line segment $A_0 A_n$ in \autoref{fig:stack_of_boxes} has length less than $\min\{\frac{r}{2}, \epsilon'\}$. Let $\theta>0$ be sufficiently small, then the linear approach $[\gamma^\theta]$ starting at $A_n$ in the direction of $\pi-\theta$ lies in $B(\gamma, t, r)$. As $\theta$ decreases, $[\gamma^\theta]$ converges to a point in $F$. 
 This shows that any open set containing $F$ has non-empty intersection with any open set containing $[\gamma]$.
 \end{proof}
\end{prop}

Before investigating more systematically translation surfaces with finite spaces of rotational components in \autoref{sec:finite_top}, we now look at examples of translation surfaces $(X, \mathcal{A})$ with singularities~$\sigma$ for which $\LT(\sigma)$ is infinite.
The first is from \cite{chamanara2004affine} and is often called the baker's map surface. 
The number of rotational components of finite and of infinite length was described in \cite{bowman_valdez2013wild}. We add a description of the topology of the space $\LT(\sigma)$ to that discussion.

\begin{exa}[Chamanara surface]\label{exa:chamanara_surface}
 Consider a square where the sides are split in segments and glued crosswise: The right half of the top is glued to the left half of the bottom, the right half of the remaining part of the top is again glued to the left half of the remaining part of the bottom, and so on as in \autoref{fig:chamanara_surface}. The same kind of gluing is also performed at the right side and the left side of the square.

 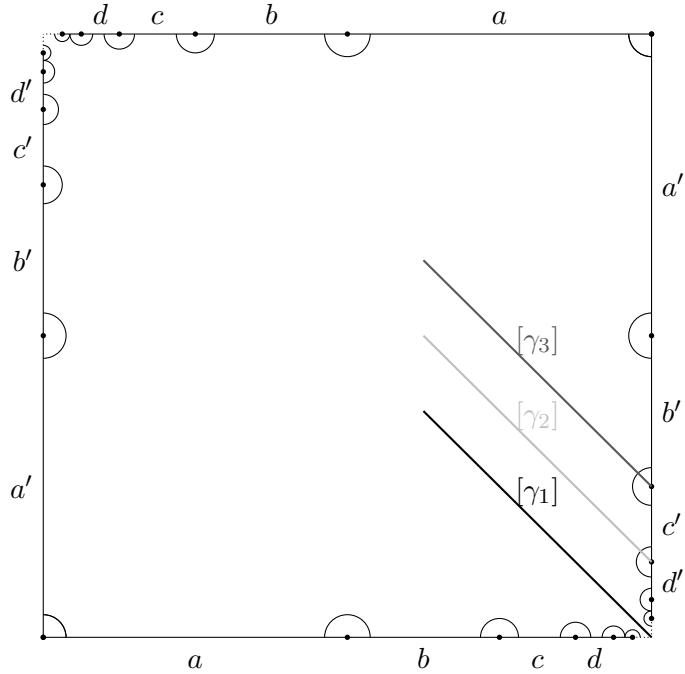
\begin{figure}
 \begin{center}
 \begin{tikzpicture}[x=1cm,y=1cm]
  \newcommand\chamanaraseitems{
   \draw (0,0) -- (7.75,0);
   \draw[densely dotted] (7.75,0) -- (8,0);
   \draw (0,0.3) arc (90:0:0.3);
   \fill (0,0) circle (1pt);
   \draw (3.7,0) arc (180:0:0.3);
   \fill (4,0) circle (1pt);
   \draw (5.75,0) arc (180:0:0.25);
   \fill (6,0) circle (1pt);
   \draw (6.8,0) arc (180:0:0.2);
   \fill (7,0) circle (1pt);
   \draw (7.35,0) arc (180:0:0.15);
   \fill (7.5,0) circle (1pt);
   \draw (7.65,0) arc (180:0:0.1);
   \fill (7.75,0) circle (1pt);
  }
  \chamanaraseitems
  \path (0,0) -- node[below=0.11cm] {$a$} (4,0) -- node[below] {$b$} (6,0) -- node[below=0.11cm] {$c$} (7,0) -- node[below] {$d$} (7.5,0);

  \begin{scope}[xscale=-1,yscale=-1,xshift=-8cm,yshift=-8cm]
   \chamanaraseitems
   \path (0,0) -- node[above] {$a$} (4,0) -- node[above] {$b$} (6,0) -- node[above] {$c$} (7,0) -- node[above] {$d$} (7.5,0);
  \end{scope}

  \begin{scope}[xscale=-1,rotate=90]
   \chamanaraseitems
   \path (0,0) -- node[left] {$a'$} (4,0) -- node[left] {$b'$} (6,0) -- node[left] {$c'$} (7,0) -- node[left] {$d'$} (7.5,0);
  \end{scope}

  \begin{scope}[xscale=-1,rotate=-90,xshift=-8cm,yshift=-8cm]
   \chamanaraseitems
   \path (0,0) -- node[right] {$a'$} (4,0) -- node[right] {$b'$} (6,0) -- node[right] {$c'$} (7,0) -- node[right] {$d'$} (7.5,0);
  \end{scope}
  
  \draw[color=black, thick] (5,3) -- node[above=0.1cm] {$[\gamma_1]$} (8,0);
  \draw[color=gray!50, thick] (5,4) -- node[above=0.1cm] {$[\gamma_2]$} (8,1);
  \draw[color=gray!70!black, thick] (5,5) -- node[above=0.1cm] {$[\gamma_3]$} (8,2);
 \end{tikzpicture}
 \label{fig:chamanara_surface}
 \caption{Chamanara surface: segments with the same letters are glued. The linear approach $[\gamma_1]$ belongs to a rotational component of length $\pi/2$, while $[\gamma_2]$ and $[\gamma_3]$ each belong to one of the two infinite-length rotational components.}
 \end{center}
 \end{figure}

 The Chamanara surface has exactly one singularity.
 From the gluings we can see that every second split point is identified, so it admits at most two singularities.
 Now the distance of these two points in the metric completion of the surface is $0$, so they are equal.
 The unique singularity $\sigma$ has two rotational components that are isometric to $\RR$ and infinitely many of finite length (see \cite{bowman_valdez2013wild}). In \autoref{fig:chamanara_surface}, representatives for three of these rotational components are shown. One more rotational component is easy to see in the figure whereas the others are easier to describe in a systematic way:
 All of the finite-length rotational components are images of each other by affine orientation-preserving homeomorphisms on the surface whose existence is guaranteed by the study in \cite{chamanara2004affine}. Furthermore, no two linear approaches in different finite-length rotational components have the same direction in $S^1=\RR/2\pi\ZZ$.
 This and the discreteness of the Veech group (see \cite[Theorem 3]{chamanara2004affine}) imply by \autoref{lem:dir_continuous} that the subspace of $\LT(\sigma)$ of all finite-length rotational components is endowed with the discrete topology.
 An argument similar to the one in \autoref{exa:geometric_series_decoration} shows that any open set containing a finite-length rotational component contains both of the infinite-length rotational components as well. For instance, the linear approaches $[\gamma_2]$ and $[\gamma_3]$ are close to $[\gamma_1]$ in \autoref{fig:chamanara_surface}.
 Moreover, each of the singletons consisting of an infinite-length rotational component is open. That can be proven directly from the definition of the topology on $\LT(\sigma)$.

 We just proved the following characterization for the open sets of $\LT(\sigma)$:
 A non-trivial set in $\LT(\sigma)$ is open if and only if it is a singleton containing an infinite-length rotational component or if it contains both infinite-length rotational components.
\end{exa}

Note that the two rotational components of length $\frac{\pi}{2}$ are readily visible in the representation of the Chamanara surface with a square, as in \autoref{fig:chamanara_surface}.
The other finite-length rotational components are harder to see in the figure, and in fact any linear approach in a rotational component of length less than $\frac{\pi}{2}$ intersects the sides of the square infinitely often. 
As we will see in the following examples, avoiding this type of behaviour simplifies the analysis of the space of rotational components. Therefore, we introduce the notion of a good cellulation.

\begin{definition}[Good cellulation]
 Let $(X, \mathcal{A})$ be a translation surface and $S$ a union of saddle connections. 
 We say $S$ is a \emph{good cellulation} if any geodesic $\gamma$ which is not in $S$ satisfies that $\{t \in \mathbb{R} :\gamma(t)\in S\}$ is a discrete set of the real line.
\end{definition} 

To be of interest, a good cellulation $S$ of a translation surface $(X, \mathcal{A})$ should also satisfy that each component of $X \setminus S$ is a piece of the Euclidean plane or of a cylinder.

As we already noted, the sides of the square in \autoref{fig:chamanara_surface} are not a good cellulation.

\begin{exa}[Star decoration]\label{exa:star_decoration}
 Let $(X,\mathcal{A})$ be a translation surface and $x\in X$. By a \emph{star decoration} of $X$ in $x$ we mean the result of cutting open all line segments starting at $x$, of length $2^{-n}$ and in the direction $\frac{m\pi}{2^n}$ such that $n\geq 1$ and $m$ is odd when $n>1$ (see \autoref{fig:star_decoration}).
 We will describe two ways to glue the branches of the star (the second will involve additional cylinders) which result in translation surfaces with non-homeomorphic spaces of rotational components.

 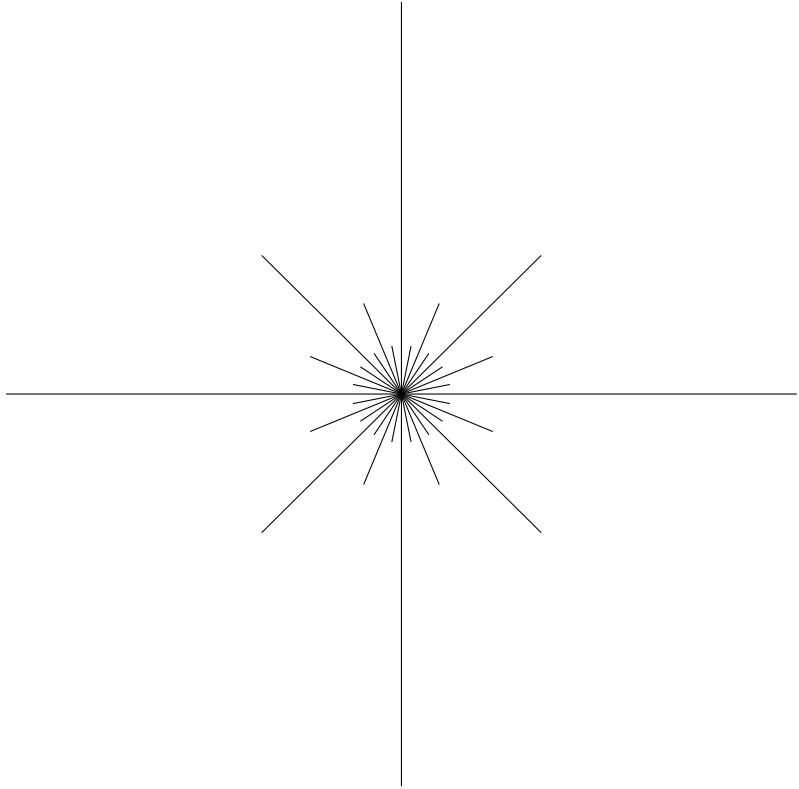
\begin{figure}[bth]
 \begin{center}
 \begin{tikzpicture}[scale=0.65]
  \draw (0,0) -- (0,8);
  \draw (0,0) -- (8,0);
  \draw (0,0) -- (-8,0);
  \draw (0,0) -- (0,-8);

  \draw (0,0) -- (45:4);
  \draw (0,0) -- (135:4);
  \draw (0,0) -- (-135:4);
  \draw (0,0) -- (-45:4);

  \draw (0,0) -- (22.5:2);
  \draw (0,0) -- (67.5:2);
  \draw (0,0) -- (112.5:2);
  \draw (0,0) -- (157.5:2);
  \draw (0,0) -- (-157.5:2);
  \draw (0,0) -- (-112.5:2);
  \draw (0,0) -- (-67.5:2);
  \draw (0,0) -- (-22.5:2);

  \draw (0,0) -- (11.25:1);
  \draw (0,0) -- (33.75:1);
  \draw (0,0) -- (56.25:1);
  \draw (0,0) -- (78.75:1);
  \draw (0,0) -- (101.25:1);
  \draw (0,0) -- (123.75:1);
  \draw (0,0) -- (146.25:1);
  \draw (0,0) -- (168.75:1);
  \draw (0,0) -- (-168.75:1);
  \draw (0,0) -- (-146.25:1);
  \draw (0,0) -- (-123.75:1);
  \draw (0,0) -- (-101.25:1);
  \draw (0,0) -- (-78.75:1);
  \draw (0,0) -- (-56.25:1);
  \draw (0,0) -- (-33.75:1);
  \draw (0,0) -- (-11.25:1);
 \end{tikzpicture}
 \label{fig:star_decoration}
 \caption{Star decoration without indications on the gluing.}
 \end{center}
 \end{figure}

 \begin{enumerate}
  \item The first kind of gluing we consider is the following:
  For each branch, glue the right side to the left side of the antipodal branch and vice versa.
  Then all tips of the branches will be identified with the center and there is only one singularity which is wild.
  For every branch, we have a rotational component which is isometric to $[0,2\pi]$. Also, for every \emph{non-dyadic} direction, i.e.\ for every direction not of the form~$\frac{m\pi}{2^n}$ for some $m, n\in\mathbb{Z}$, we have a linear approach starting in the center of the star.
  The corresponding rotational component consists only of this point.

  The problem with this gluing is that the union of all the branches is not a good cellulation. 
  In fact, there are additional rotational components that are not easy to see in the figure. 
  It is possible to find such a rotational component by starting with any point not on a branch, and considering geodesics passing through that point.
  One can then inductively find nested open sets $U_n$ of directions for which the geodesics intersect $n$ smaller and smaller branches.
  By construction the intersection of these nested open intervals is a singleton, containing a direction for which the geodesic reaches the center of the star in finite time (i.e.\ defines a linear approach) while passing through infinitely many branches (see \autoref{fig:star_decoration1}).
  
   \begin{figure}[bthp]
   	\begin{center}
   		\begin{tikzpicture}[scale=0.8]
   		\draw[color=gray] (0,0) -- (0,8);
   		\draw[color=gray, name path=branch0] (0,0) -- (8,0);
   		\draw[color=gray] (0,0) -- (-8,0);
   		\draw[color=gray] (0,0) -- (0,-8);
   		
   		\draw[color=gray] (0,0) -- (45:4);
   		\draw[color=gray] (0,0) -- (135:4);
   		\draw[color=gray, name path=branch1] (0,0) -- (-135:4);
   		\draw[color=gray] (0,0) -- (-45:4);
   		
   		\draw[color=gray, very thin, name path=branch2] (0,0) -- (22.5:2);
   		\draw[color=gray, very thin] (0,0) -- (67.5:2);
   		\draw[color=gray, very thin] (0,0) -- (112.5:2);
   		\draw[color=gray, very thin] (0,0) -- (157.5:2);
   		\draw[color=gray, very thin] (0,0) -- (-157.5:2);
   		\draw[color=gray, very thin] (0,0) -- (-112.5:2);
   		\draw[color=gray, very thin] (0,0) -- (-67.5:2);
   		\draw[color=gray, very thin] (0,0) -- (-22.5:2);
   		
   		\draw[color=gray, ultra thin] (0,0) -- (11.25:1);
   		\draw[color=gray, ultra thin] (0,0) -- (33.75:1);
   		\draw[color=gray, ultra thin] (0,0) -- (56.25:1);
   		\draw[color=gray, ultra thin] (0,0) -- (78.75:1);
   		\draw[color=gray, ultra thin] (0,0) -- (101.25:1);
   		\draw[color=gray, ultra thin] (0,0) -- (123.75:1);
   		\draw[color=gray, ultra thin] (0,0) -- (146.25:1);
   		\draw[color=gray, ultra thin] (0,0) -- (168.75:1);
   		\draw[color=gray, ultra thin] (0,0) -- (-168.75:1);
   		\draw[color=gray, ultra thin, name path=branch3] (0,0) -- (-146.25:1);
   		\draw[color=gray, ultra thin] (0,0) -- (-123.75:1);
   		\draw[color=gray, ultra thin] (0,0) -- (-101.25:1);
   		\draw[color=gray, ultra thin] (0,0) -- (-78.75:1);
   		\draw[color=gray, ultra thin] (0,0) -- (-56.25:1);
   		\draw[color=gray, ultra thin] (0,0) -- (-33.75:1);
   		\draw[color=gray, ultra thin] (0,0) -- (-11.25:1);
   		
   		\path (4,1.7) coordinate (q0);
   		\path[name path=path0] (q0) -- ++(-60:4);
   		\draw[thick, name intersections={of=branch0 and path0, by={p0}}] (q0) -- (p0);
   		\path (p0) ++ (-180:8) coordinate (q1);
   		\path[name path=path1] (q1) -- ++(-60:4);
   		\draw[thick, name intersections={of=branch1 and path1, by={p1}}] (q1) -- (p1);
   		\path (p1) ++(45:4) coordinate (q2);
   		\path[name path=path2] (q2) -- ++(-60:4);
   		\draw[thick, name intersections={of=branch2 and path2, by={p2}}] (q2) -- (p2);
   		\path (p2) ++ (-157.5:2) coordinate (q3);
   		\path[name path=path3] (q3) -- ++(-60:4);
   		\draw[thick, name intersections={of=branch3 and path3, by={p3}}] (q3) -- (p3);
   		\path (p3) ++ (33.75:1) coordinate (q4);
   		\draw[thick] (q4) -- ++ (-60:0.03);
   		\end{tikzpicture}
   		\label{fig:star_decoration1}
   		\caption{In the star decoration without additional cylinders, the geodesic curve starts in the center of the star and intersects infinitely many branches.}
   	\end{center}
   \end{figure}
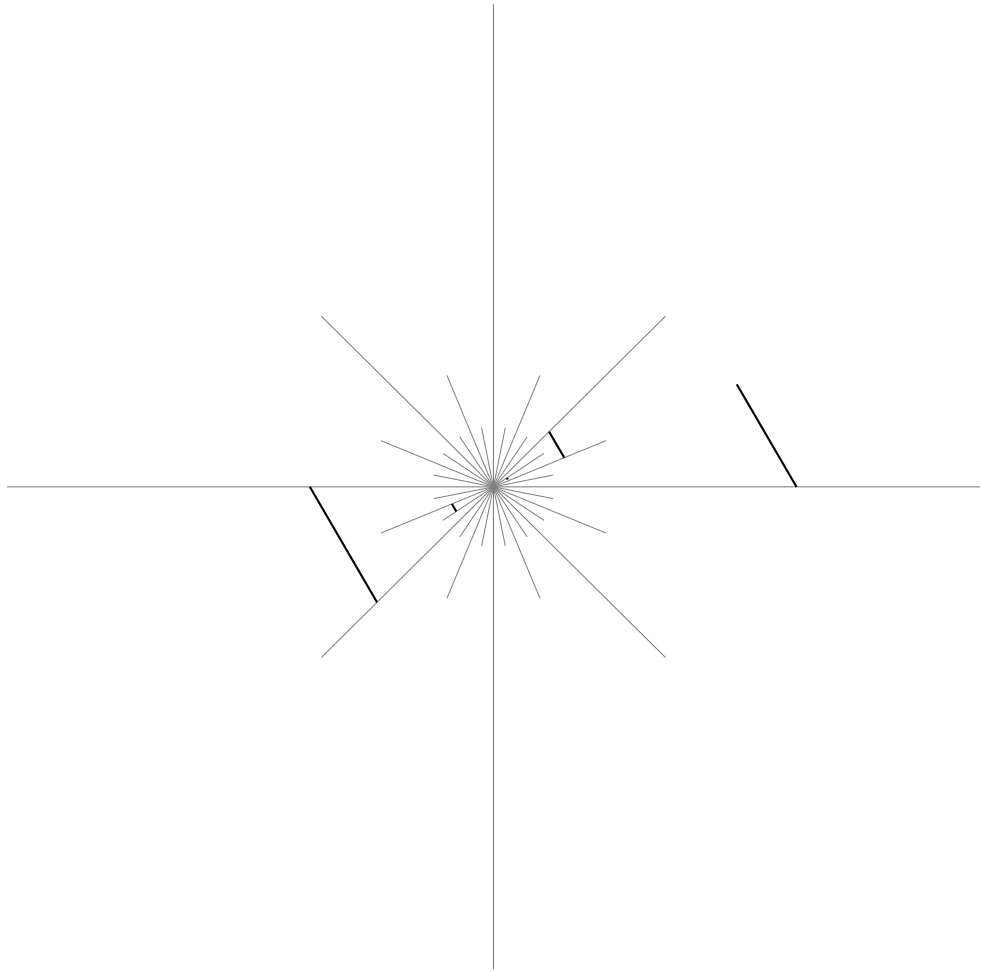

  \item We now want to turn the surface with the star decoration into a surface with a good cellulation.
  For every branch, consider a rectangle as in \autoref{fig:neighbored_squares_glued}, rotate it so that the former vertical sides are parallel to the branch and glue these sides to the sides of the branch. This rectangle consists of $2^n + 1$ squares of side length $2^{-n}$, with the tops and bottoms of the squares glued crosswise, except for the middle square where top and bottom are identified.

  \begin{figure}
  \begin{center}
  	\newlength\zzz
  	\settoheight\zzz{$b$}
  \begin{tikzpicture}[scale=1.5]
   \draw (0,0) -- node[below, text height=\zzz]{$b$} 
   (1,0) -- node[below, text height=\zzz]  {$a$}
    (2,0) -- node[below, text height=\zzz] {$d$}
    (3,0) -- node[below, text height=\zzz] {$c$}
    (4,0) -- node[below, text height=\zzz] {$e$}
    (5,0) -- node[below, text height=\zzz] {$g$}
    (6,0) -- node[below, text height=\zzz] {$f$}
    (7,0) -- node[below, text height=\zzz] {$i$}
    (8,0) -- node[below, text height=\zzz] {$h$}
    (9,0) --
    (9,1) -- node[above] {$i$}
    (8,1) -- node[above] {$h$}
    (7,1) -- node[above] {$g$}
    (6,1) -- node[above] {$f$}
    (5,1) -- node[above] {$e$}
    (4,1) -- node[above] {$d$}
    (3,1) -- node[above] {$c$}
    (2,1) -- node[above] {$b$}
    (1,1) -- node[above] {$a$}
    (0,1) -- (0,0);
   \draw (1,0) -- (1,1);
   \draw (2,0) -- (2,1);
   \draw (3,0) -- (3,1);
   \draw (4,0) -- (4,1);
   \draw (5,0) -- (5,1);
   \draw (6,0) -- (6,1);
   \draw (7,0) -- (7,1);
   \draw (8,0) -- (8,1);
  \end{tikzpicture}
  \label{fig:neighbored_squares_glued}
  \caption{Segments with the same letters are glued.}
  \end{center}
  \end{figure}
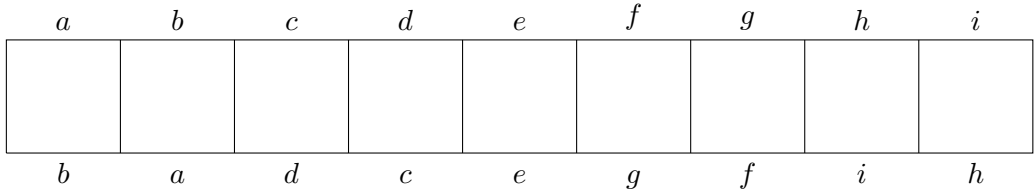

  Now the metric completion of the resulting surface is compact if the star decoration is performed on a compact surface. This is because any $\epsilon$--neighborhood of the center of the star must contain all but finitely many branches and the rectangles glued to them. 
  Let $S$ be the union of all the branches of the star.
  By construction, $S$ is a good cellulation.

  As in the previous case, we have exactly one singularity $\sigma$ as all tips of the branches are identified with the center as well as the vertices of the squares in the additional rectangles.
 \end{enumerate}
\end{exa}

Whereas for the star decoration without additional cylinders, the space of rotational components is very complicated to describe, it is more feasible for the star decoration with additional cylinders. In fact, the space of rotational components is more handsome than the construction suggests.

\begin{prop}[Space of rotational components for star decoration]
 For a translation surface with a star decoration with additional cylinders and singularity~$\sigma$, the space $\LT(\sigma)$ of rotational components is homeomorphic to $S^1$.

\begin{proof}
 Because $S$ is a good cellulation, we see that there are exactly four types of linear approaches:
  \begin{enumerate}[a)]
   \item rays $\gamma$ starting from the center of the star going in a non-dyadic direction, and such that $\gamma(t)$ is not in the interior of a rectangle for small $t$,
   \item rays $\gamma$ starting from the center of the star going in a dyadic direction, and such that~$\gamma(t)$ is not in the interior of a rectangle for small $t$,
   \item rays $\gamma$ starting from the tip of a branch of the star, and such that $\gamma(t)$ is not in the interior of a rectangle for small $t$,
   \item rays $\gamma$ starting at a vertex of a square in a rectangle and such that $\gamma(t)$ stays in the interior of this rectangle for small $t$.
  \end{enumerate}

  Now there are two types of rotational components:
  \begin{enumerate}[A)]
   \item singletons consisting of a linear approach of type a), 
   \item finite-length rotational components isometric to $[0, (2^n+2) \cdot 2\pi]$ via a map $\phi$ such that $\phi(0)$ and $\phi( (2^n+2) \cdot 2\pi)$ are of type b) and go in the same (dyadic) direction, $\phi(x)$ is of type c) for $x \in [(2^n+1) \cdot \pi, (2^n+3) \cdot \pi]$ and of type d) otherwise.
  \end{enumerate}
  Therefore, we can define a one-to-one map $P$ from $\LT(\sigma)$ to $S^1$ that associates to a rotational component the direction of its element of type a) or b).
  We will now show that $P$ is a homeomorphism by proving that it is open and continuous.

  To prove that $P$ is open, we show that
  for every rotational component $c$ (seen as a class of linear approaches in $\LC(\sigma)$), the image of any neighborhood of $c$ in $\LC(\sigma)$ is a neighborhood of $P(c)$ in $S^1$.
  Assume first that $c$ is a singleton consisting of a linear approach $[\gamma]$ of type a).
  Any neighborhood of $[\gamma]$ contains some open set $B(\gamma, t, r)$ as we have shown in \autoref{lem:subbasis_LCpoint}. 
  Let $\theta$ be the direction of $[\gamma]$, that is $\gamma(t) = te^{\mathrm{i}\theta}$.
  Choose $\epsilon$ and $\eta$ small enough such that
  \begin{equation*}
   \{(t+\rho)e^{\mathrm{i}(\theta+\varphi)} : \rho \in [0, \eta], \varphi\in (-\epsilon,\epsilon)\}
  \end{equation*}
  is contained in the disk with center $\gamma(t)$ and radius $r$.
  By decreasing $\epsilon$ if necessary, we can assume that no branch with direction in $(\theta-\epsilon,\theta+\epsilon)$ is longer than $\eta$. 
  Then any ray of type a) with direction in $(\theta-\epsilon,\theta+\epsilon)$ is in $B(\gamma, t, r)$, as is any ray of type c) in direction $\varphi\in(\theta-\epsilon,\theta+\epsilon)$ starting on a branch with the same direction $\varphi$.
  Thus $P(B(\gamma, t, r))$ contains $(\theta-\epsilon,\theta+\epsilon)$ so it is a neighborhood of $\theta = P(c)$.
  If $P(c)$ is a dyadic direction, we can conclude with a similar argument for the two linear approaches of type~b) in~$c$.

  To prove that $P$ is continuous, we show that its pre-composition with the quotient map $\LC(\sigma) \rightarrow \LT(\sigma)$ is continuous. Denote the composition as $P' \colon \LC(\sigma) \to S^1$.
  Let $[\gamma] \in \LC(\sigma)$ and $\gamma$ a representative of $[\gamma]$. We have to show that there exist $\epsilon > 0$, $t > 0$, and $r > 0$ so that for any $[\gamma']\in B(\gamma, t, r)$, the image $P'([\gamma'])$ is contained in $(P'([\gamma])-\epsilon,P'([\gamma])+\epsilon)$. As this is done by a case-by-case analysis depending on the type of both $[\gamma]$ and $[\gamma']$, we will not carry this out completely but give instructions how to choose $\epsilon$, $t$, and $r$.
  
  For the fixed $[\gamma] \in \LC(\sigma)$, choose $\epsilon > 0$ but very close to $0$.
  Furthermore, choose $t > 0$ small enough so that there is a branch of length greater than $10t$ in a direction in $(P'([\gamma])-\epsilon,P'([\gamma]))$, and so that there is a branch of length greater than $10t$ in a direction in $(P'([\gamma]),P'([\gamma])+\epsilon)$. This choice of $t$ will prevent the images of linear approaches of type c) to be in $(P'([\gamma])-\epsilon,P'([\gamma]) + \epsilon)$ when the corresponding branch on which tip the linear approach is starting does not already have a direction in $(P'([\gamma])-\epsilon,P'([\gamma]) + \epsilon)$.
  Moreover, let $r \in (0,\frac{1}{10}t)$ be small enough such that $B(\gamma(t),r)$ does not intersect any branch except possibly the one $\gamma$ lies on.
  By these conditions it follows that for any $[\gamma']\in B(\gamma, t, r)$, the image $P'([\gamma'])$ is contained in $(P'([\gamma])-\epsilon,P'([\gamma])+\epsilon)$.
 
  For example, suppose both $[\gamma]$ and $[\gamma']$ are of type a) but $P'([\gamma'])$ is not contained in $(P'([\gamma])-\epsilon,P'([\gamma])+\epsilon)$. Then by the choice of $t$ there must be a branch of length more than $10t$ which intersects the geodesic from $\gamma(t)$ to $\gamma'(t)$. Hence the distance between the points $\gamma(t)$ and $\gamma'(t)$ must be greater than $r$, which is a contradiction.
  Therefore, $P$ is a homeomorphism as was to be shown.
\end{proof}
\end{prop}

A modification of the previous example will be very helpful in \autoref{sec:finite_top} so we introduce it here.
It will serve as a building block in the proof of \autoref{thm:finite_top}.

\begin{exa}[Shrinking star decoration]\label{exa:shrinking_star_decoration}
 Consider $\mathbb{R}^2$ with geodesic segments $l_n$ from the origin $(0,0)$ to $(2^{-n}\sin(\frac{\pi}{2^n}), 2^{-n}\cos(\frac{\pi}{2^n}))$ for $n\geq 1$ as in \autoref{fig:shrinking_star_decoration}.
 Make a slit along each $l_n$, glue the two sides to the two smaller sides of a $(1+2^{-n})$--by--$2^{-n}$ rectangle as in \autoref{fig:neighbored_squares_glued} and then glue the other sides of the squares in this rectangle in a crosswise way as before.

 \begin{figure}[]
 \begin{center}
 \begin{tikzpicture}[scale=2]
  \draw (0:0) -- (0:4);
  \draw (0:0) -- (45:2);
  \draw (0:0) -- (67.5:1);
  \draw (0:0) -- (78.75:0.5);
  \draw (0:0) -- (84.375:0.25);
  \draw (0:0) -- (87.1875:0.125);
 \end{tikzpicture}
 \label{fig:shrinking_star_decoration}
 \caption{Shrinking star decoration.}
 \end{center}
 \end{figure}
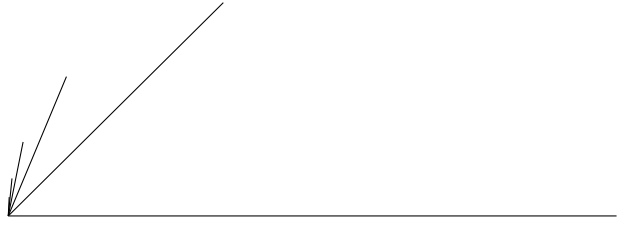

 In this case, we have one wild singularity with one rotational component, whose translation structure is isomorphic to $[0,\infty)$.
\end{exa}

\section[Relations between the translation structure on rotational components and the topology on \texorpdfstring{$\LT(\sigma)$}{the space of rotational components}]{Relations between the translation structure on \linebreak rotational components and the topology on \texorpdfstring{$\LT(\sigma)$}{the space of rotational components}}

In this section we describe relations between two topologies that are associated to the rotational components of a singularity $\sigma$. These are 
\begin{enumerate}
 \item the translation structure on rotational components, seen as classes of linear approaches, and
 \item the topology on $\LT(\sigma)$, coming from the uniform metric on the spaces $\LC^\epsilon(X)$.
\end{enumerate}

We start with a topological property of a rotational component in $\LT(\sigma)$ that can be concluded from a property of the translation structure on this rotational component.

\begin{prop}[Rotational components isometric to $\RR$ are open points]\label{prop:biinfinite_means_open}
 Let $(X,\mathcal{A})$ be a translation surface with discrete singularities and $\sigma$ be a wild singularity of $(X,\mathcal{A})$. If $c$ is a rotational component of $\sigma$ which is isometric to $\RR$ then $\{c\}$ is open in $\LT(\sigma)$.

\begin{proof}
 Let $c \in \LT(\sigma)$ be a rotational component whose one-dimensional translation structure is isometric to $\RR$.
 The set $\{c\}$ is open in $\LT(\sigma)$ if the set of linear approaches contained in $c$ is open in $\LC(\sigma)$. 
 To prove the statement of the proposition, we give an open neighborhood for every linear approach contained in $c$ that does not contain linear approaches from other rotational components. The union of these neighborhoods is then an open set in $\LC(\sigma)$.

 Let $[\gamma]$ be a linear approach in $c$ with $\gamma \in \LC^\epsilon(\sigma)$ a representative for a suitable $\epsilon>0$. Because $c$ is isometric to $\RR$, there exists an angular sector $([-2\pi,2\pi],\ln\epsilon,i_{\ln\epsilon})$ such that
 \begin{equation*}
  i_{\ln\epsilon} ( (-\infty, \ln\epsilon) \times \{ 0 \} ) = \gamma((0,\epsilon)) .
 \end{equation*}
 By \autoref{lem:subbasis_LCpoint}, $B(\gamma, \frac{\epsilon}{4}, \frac{\epsilon}{4})$ is an open neighborhood of $[\gamma]$ in $\LC(\sigma)$.
 Recall that it is the set of all linear approaches with a representative $\gamma'$ for which $\gamma'(\frac{\epsilon}{4})$ is contained in the disk around $\gamma(\frac{\epsilon}{4})$ of radius $\frac{\epsilon}{4}$. We will show that the linear approaches in $B(\gamma, \frac{\epsilon}{4}, \frac{\epsilon}{4})$ are also contained in $c$.

 The set
 \begin{equation*}
  i_{\ln\epsilon}((-\infty,\ln\epsilon)\times(-\pi,\pi))
 \end{equation*}
 in $X$ is isometric to an open $\epsilon$--disk in $\mathbb{R}^2$ with a slit removed (see  \autoref{fig:isometric_to_RR_means_open_proof_epsilon_disk}).

 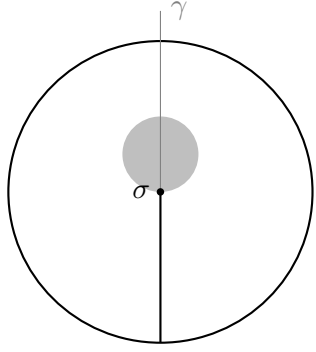
\begin{figure}[]
 \begin{center}
 \begin{tikzpicture}[scale=2]
  \draw[thick](0,0) circle(1);
  \fill[color=gray!50](0,0.25) circle(0.25);
  \draw[thick] (0,0) -- (0,-1);
  \draw[color=gray] (0,0)--(0,1.2) node[right]{$\gamma$};
  \fill (0,0) node[left] {$\sigma$} circle (0.7pt);
 \end{tikzpicture}
 \label{fig:isometric_to_RR_means_open_proof_epsilon_disk}
 \caption{The set $i_{\ln\epsilon}((-\infty,\ln\epsilon)\times(-\pi,\pi))$ with $\gamma$ and the shaded region $B(\gamma(\frac{\epsilon}{4}),\frac{\epsilon}{4})$.}
 \end{center}
 \end{figure}

 Let $s \colon (0, \frac{\epsilon}{4}] \to X$ be a geodesic curve with $s(\frac{\epsilon}{4}) \in B(\gamma(\frac{\epsilon}{4}),\frac{\epsilon}{4})$. Then the image of $s$ either intersects the slit or it is contained in $i_{\ln\epsilon}( (-\infty,\ln\epsilon)\times(-\pi,\pi) )$.

 In the first case, it can hit the slit only at the dashed segments
 \begin{equation*}
  i_{\ln\epsilon}((-\infty,\ln\frac{\epsilon}{4})\times\{-\pi\})
  \quad \text{or} \quad
  i_{\ln\epsilon}((-\infty,\ln\frac{\epsilon}{4})\times\{\pi\})
 \end{equation*}
 in \autoref{fig:isometric_to_RR_means_open_proof_half_disks}.

 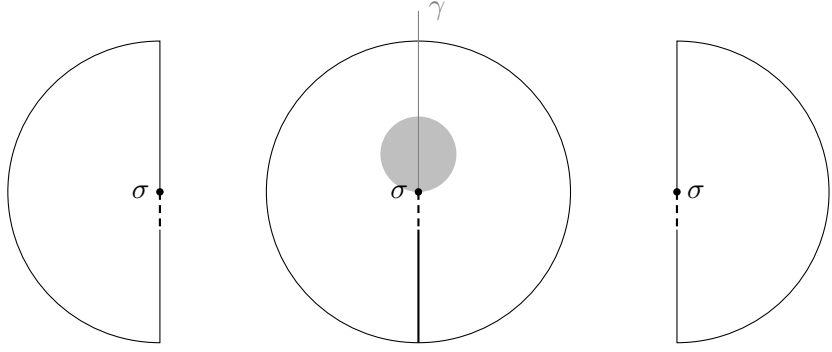
\begin{figure}
 \begin{center}
 \begin{tikzpicture}[scale=2]
 \draw (-1.7,0) -- (-1.7,1) arc (90:270:1) -- (-1.7,-0.25);
  \draw[densely dashed,thick](-1.7,0)--(-1.7,-0.25);
  \fill (-1.7,0) node[left] {$\sigma$} circle (0.7pt);

  \draw (0,0) circle(1);
  \fill[color=gray!50](0,0.25) circle(0.25);
  \draw[thick](0,-0.25)--(0,-1);
  \draw[densely dashed, thick](0,0)--(0,-0.25);
  \draw[color=gray] (0,0)--(0,1.2) node[right]{$\gamma$};
  \fill (0,0) node[left] {$\sigma$} circle (0.7pt);
  
  \draw (1.7,0) -- (1.7,1) arc (90:-90:1) -- (1.7,-0.25);
  \draw[densely dashed, thick] (1.7,0)--(1.7,-0.25);
  \fill (1.7,0) node[right] {$\sigma$} circle (0.7pt);
 \end{tikzpicture}
 \label{fig:isometric_to_RR_means_open_proof_half_disks}
 \caption{From left to right, we have the three sets $i_{\ln\epsilon}((-\infty, \ln\epsilon)\times(-2\pi,-\pi))$, $i_{\ln\epsilon}((-\infty, \ln\epsilon)\times(-\pi,\pi))$ and $i_{\ln\epsilon}((-\infty,\ln\epsilon)\times(\pi,2\pi))$.}
 \end{center}
 \end{figure}

 If it hits $i_{\ln\epsilon}((-\infty,\ln\frac{\epsilon}{4})\times\{-\pi\})$ or $i_{\ln\epsilon}((-\infty,\ln\frac{\epsilon}{4})\times\{\pi\})$, consider
 \begin{equation*}
  i_{\ln\epsilon}((-\infty, \ln\epsilon)\times(-2\pi,-\pi))
  \quad \text{or} \quad
  i_{\ln\epsilon}((-\infty,\ln\epsilon) \times (\pi,2\pi)),
 \end{equation*}
 which are open half disks with radius $\epsilon$ as in \autoref{fig:isometric_to_RR_means_open_proof_half_disks}.
 
 A geodesic segment shorter than $\frac{\epsilon}{4}$ starting on either of the dashed segments, heading inwards, has its endpoint in the interior of this half disk,
 hence $s$ cannot start in a singularity. This means that $s$ is not contained in $\LC^{\frac{\epsilon}{4}} (\sigma)$ and does not define a linear approach to $\sigma$.

 In the second case, if the geodesic segment $s$ stays in $i_{\ln\epsilon}( (-\infty,\ln\epsilon)\times(-\pi,\pi) )$ and cannot be extended, it can only define a linear approach of $\sigma$ if it starts at the center. Hence, it must represent a linear approach in the same rotational component $c$ as $[\gamma]$.
\end{proof}
\end{prop}

Note that in the proof of the previous proposition, the condition on the rotational component being isometric to $\mathbb{R}$ was only used to ensure that for each linear approach there exists an angular sector $(I, c, i_c)$ with sufficiently large interval $I$. We extract this as a separate lemma on linear approaches as we will use it in the proof of \autoref{thm:info_loss}.

\begin{lem}[Neighborhoods of linear approaches with distance greater $\pi$ to the boundary] \label{lem:neighborhood_linear_approaches}
 Let $(X,\mathcal{A})$ be a translation surface with discrete singularities, $\sigma$ be a wild singularity of $(X,\mathcal{A})$, and $[\gamma] \in \LC(\sigma)$ a linear approach that has distance greater than $\pi$ to the boundary of the rotational component $\overline{[\gamma]}$ with the translation structure.
 
 Then $[\gamma]$ has an open, path-connected neighborhood in $\LC(\sigma)$ which is entirely contained in $\overline{[\gamma]}$.
 
\begin{proof}
 Let $\gamma \in \LC^\epsilon(\sigma)$ be a representative of $[\gamma]$ for a suitable $\epsilon>0$. As the distance of $[\gamma]$ to the boundary is greater than $\pi$, there exists a $\delta > 0$ and an angular sector $([-\pi-\delta,\pi+\delta],\ln\epsilon,i_{\ln\epsilon})$ such that
 \begin{equation*}
  i_{\ln\epsilon} ( (-\infty, \ln\epsilon) \times \{ 0 \} ) = \gamma((0,\epsilon)) .
 \end{equation*}
 We follow the same proof as in \autoref{prop:biinfinite_means_open} but replace the open half disks in \autoref{fig:isometric_to_RR_means_open_proof_half_disks} by sectors of central angle $\delta$. Furthermore, we replace $B(\gamma, \frac{\epsilon}{4}, \frac{\epsilon}{4})$ by $B(\gamma, \frac{\epsilon}{4}, \epsilon')$ with $\epsilon' \in (0,\frac{\epsilon}{4})$ small enough so that $B(\gamma(\frac{\epsilon}{4}), \epsilon')$ is contained in $i_{\ln\epsilon} ( (-\infty, \ln\epsilon) \times (-\delta, \delta) )$. This implies that no geodesic $s$ as in the proof of \autoref{prop:biinfinite_means_open} can start outside of the sectors (see \autoref{fig:neighborhood_linear_approaches_proof}).
 
 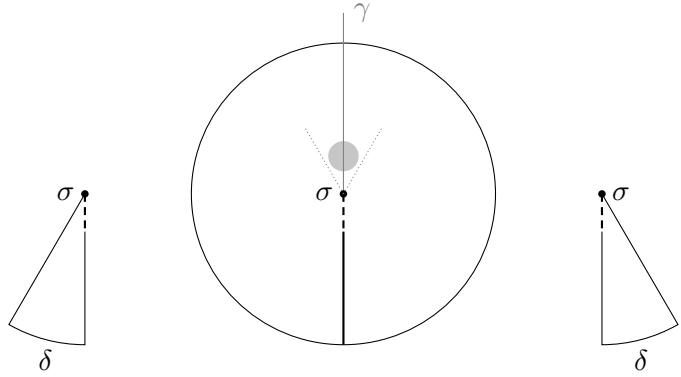
\begin{figure}
 \begin{center}
 \begin{tikzpicture}[scale=2]
 \draw (-1.7,0) -- ++(240:1) arc (240:270:1) node[midway, below] {$\delta$} -- (-1.7,-0.25);
  \draw[densely dashed,thick](-1.7,0)--(-1.7,-0.25);
  \fill (-1.7,0) node[left] {$\sigma$} circle (0.7pt);

  \draw (0,0) circle(1);
  \fill[color=gray!45](0,0.25) circle(0.1);
  \draw[thick](0,-0.25)--(0,-1);
  \draw[densely dashed, thick](0,0)--(0,-0.25);
  \draw[color=gray] (0,0)--(0,1.2) node[right]{$\gamma$};
  \fill (0,0) node[left] {$\sigma$} circle (0.7pt);
  \draw[color=gray, densely dotted] (0,0) -- ++(60:0.5);
  \draw[color=gray, densely dotted] (0,0) -- ++(120:0.5);
  
  \draw (1.7,0) -- ++(-60:1) arc (-60:-90:1) node[midway, below] {$\delta$} -- (1.7,-0.25);
  \draw[densely dashed, thick] (1.7,0)--(1.7,-0.25);
  \fill (1.7,0) node[right] {$\sigma$} circle (0.7pt);
 \end{tikzpicture}
 \label{fig:neighborhood_linear_approaches_proof}
 \caption{The dotted line indicates that in this variation of the proof no segment~$s$ can start outside of the sectors.}
 \end{center}
 \end{figure}
 
 As before, $B(\gamma, \frac{\epsilon}{4}, \epsilon')$ is an open neighborhood of $[\gamma]$ which is entirely contained in~$\overline{[\gamma]}$. The explicit description of the linear approaches in $B(\gamma, \frac{\epsilon}{4}, \epsilon')$ shows that this neighborhood is even path-connected.
\end{proof}
\end{lem}

In general, the converse statement of \autoref{prop:biinfinite_means_open} is not true. The stack of boxes from \autoref{exa:stack_of_boxes}, for instance, has a wild singularity $\sigma$ with exactly one rotational component~$c$. Hence, $\{c\}$ is open in $\LT(\sigma) = \{c\}$. On the other hand, the rotational component~$c$ has a one-dimensional translation structure which is isometric to $(0,\infty)$.
However, a slightly weaker statement than the converse of \autoref{prop:biinfinite_means_open} is true.

\begin{prop}[Open rotational components have infinite length]\label{prop:open_means_infinite_length}
 Let $(X,\mathcal{A})$ be a translation surface with discrete singularities and $\sigma$ be a wild singularity of $(X, \mathcal{A})$. If $c$ is a rotational component such that $\{c\}$ is open in $\LT(\sigma)$ then it is isometric to an interval in $\RR$ of infinite length.

\begin{proof}
 Consider a rotational component $c$ of finite length. We show that there exists a sequence of linear approaches not contained in $c$ that converges to a linear approach contained in $c$.
 This implies that $\{c\}$ is not open in $\LT(\sigma)$.

 As $c$ is of finite length, it is isometric to an interval with endpoints $a,b\in \mathbb{R}$ (the interval can be open, closed or half-open, half-closed). Choose a linear approach $[\gamma]$ which corresponds to a point in the interval that differs at most $\frac{\pi}{4}$ from $a$.
 Choose a representative $\gamma$ of $[\gamma]$, a regular point $\gamma(t_0)$ on it,
 and an $\epsilon_1 > 0$ small enough so that $B(\gamma(t_0), \epsilon_1) \subseteq X$ is isometric to a flat disk of radius $\epsilon_1$ (see \autoref{fig:open_means_infinite_length_proof} for a sketch).

 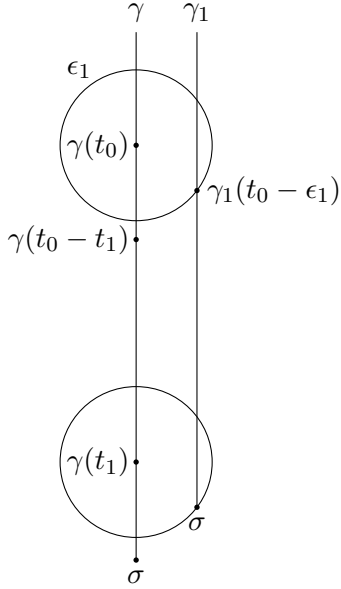
\begin{figure}
 \begin{center}
 \begin{tikzpicture}[scale=1]
  \draw (0,1.3) circle (1cm);
  \draw (0,5.5) circle (1cm);
  \fill(0,1.3) circle (1pt);
  \fill(0,5.5) circle (1pt);
  \fill(0,0) node[below] {$\sigma$} circle (1pt);
  \fill(0.8,0.7) node[below] {$\sigma$} circle (1pt);
  \fill(0,4.25) circle (1pt);
  \fill(0.8,4.9) circle (1pt);
  \draw (0,0) -- (0,1.3) node[left=-0.05cm]{$\gamma(t_1)$} -- (0,4.25) node[left=-0.05cm]{$\gamma(t_0-t_1)$} -- (0,5.5) node[left=-0.05cm]{$\gamma(t_0)$} -- (0,7) node[above]{$\gamma$};
  \draw (0.8,0.7) -- (0.8,4.9) node[right]{$\gamma_1(t_0-\epsilon_1)$} -- (0.8,7) node[above]{$\gamma_1$};
  \draw (-0.75,6.5) node {$\epsilon_1$};
 \end{tikzpicture}
 \label{fig:open_means_infinite_length_proof}
 \caption{For the linear approach $[\gamma]$ there exists a linear approach $[\gamma_1]$ so that the representatives $\gamma$ and $\gamma_1$ are parallel with distance at most $\epsilon_1$.}
 \end{center}
 \end{figure}

 Let
 \begin{equation*}
  t_1 \coloneqq \inf\{t > 0 : B(\gamma(t), \epsilon_1)\text{ is isometric to a flat disk of radius }\epsilon_1 \} .
 \end{equation*}
 By definition, we have $t_1 \leq t_0$ but we also have $t_1 \geq \epsilon_1 > 0$.
 Indeed, the disk $B(\gamma(\epsilon_1), \epsilon_1)$ cannot be isometric to a flat disk since the point corresponding to $[\gamma]$ in the interval isometric to $c$ would have distance at least $\frac{\pi}{2}$ from $a$.

 By definition of $t_1$, the boundary of the disk $B(\gamma(t_1), \epsilon_1)$ contains a singularity. As the set of singularities of $(X, \mathcal{A})$ is discrete, for $t_0$ and $\epsilon_1$ small enough it is the singularity~$\sigma$ which must be contained in the boundary of the disk $B(\gamma(t_1), \epsilon_1)$.
 This determines a geodesic curve $\gamma_1$ in $X$ starting at $\sigma$ on the boundary of $B(\gamma(t_1),\epsilon_1)$ going in the same direction as $\gamma$. By definition of $t_1$, the curve $\gamma_1$ can be defined at least for times in $(0, t_0-t_1)$ (see \autoref{fig:open_means_infinite_length_proof}).
 Therefore the distance between $\gamma$ and $\gamma_1$ in $\LC^{t_0-t_1}(\sigma)$ with the uniform metric is no more than $t_1+\epsilon_1$.

 Repeat this construction with $t_1' \leq t_1$ and $\epsilon_2 >0$ small enough so that $\gamma_1$ is not contained in the $(t_1'+\epsilon_2)$--neighborhood of $\gamma$ in $\LC^{t_0-t_1}(\sigma)$.
 We get a different linear approach $[\gamma_2]$ closer to $[\gamma]$.
 By repeating this construction, we iteratively obtain infinitely many different linear approaches $([\gamma_n])_{n \geq 1}$ which all have the same direction.
 Only finitely many of them can be contained in $c$ as $c$ is of finite length.
 Thus there is a number $N\geq 1$ for which $[\gamma_n]$ is not contained in $c$ for $n>N$, and $([\gamma_n])_{n> N}$ converges to $[\gamma]$ in $c$.
 Therefore $\{c\}$ is not open in $\LT(\sigma)$.
\end{proof}
\end{prop}

The proofs of \autoref{prop:biinfinite_means_open} and \autoref{prop:open_means_infinite_length} suggest that the interesting linear approaches -- in terms of the topology of a rotational component as subspace of $\LC(\sigma)$~-- are the linear approaches that differ at most $\pi$ from the boundary of the rotational component with the translation structure. We will use this insight in the constructions in the proofs of \autoref{thm:finite_top} and \autoref{thm:info_loss}.

The statements in the last two propositions are on rotational components of infinite length. Note that not all singularities have a rotational component of infinite length. For example the singularity of the star decoration in \autoref{exa:star_decoration} (ii) has infinitely many rotational components but all of them are of finite length. However, it is not possible that all rotational components are shorter than $\pi$ as we will see in the next proposition.

\begin{prop}[Rotational components of length no less than $\pi$ are dense]\label{prop_rotational_components_of_length_at_least_pi_are_dense}
Let $(X,\mathcal{A})$ be a translation surface and $\sigma$ a wild singularity.
Then the set of rotational components of length greater than or equal to $\pi$ is dense in $\LT(\sigma)$.

\begin{proof}
 Let $\mathcal{S}$ be the set of rotational components that contain a linear approach $\gamma$, such that for some small $t>0$, $d(\gamma(t), \sigma)=t$. With similar arguments as before, it is clear that all elements of $\mathcal{S}$ have length at least $\pi$. We show that $\mathcal{S}$ is dense in $\LT(\sigma)$.

 Let $c \in \LT(\sigma) \setminus \mathcal{S}$ be a rotational component with a linear approach $[\gamma] \in c$. Then the argument in the proof of \autoref{prop:open_means_infinite_length} shows that any neighborhood of $[\gamma]$ must contain some linear approach from a rotational component in $\mathcal{S}$. Hence $\mathcal{S}$ is dense.
\end{proof}
\end{prop}

A consequence of the fact that the set $\mathcal{S}$ from the preceding proof is dense is that the space of rotational components is separable.

\begin{cor}[$\LT(\sigma)$ is separable] \label{cor:separable}
 Let $(X,\mathcal{A})$ be a translation surface and $\sigma$ a wild singularity of $(X, \mathcal{A})$. Then the space $\LT(\sigma)$ of rotational components is separable.

\begin{proof}
 Because $X$ is separable, there exists a dense, countable subset $Y$ of $X$.
 Furthermore, for every $x \in Y$, consider the set of linear approaches $[\gamma]\in \LC(\sigma)$ with representative~$\gamma$ such that $\gamma( d(\sigma, x) ) = x$. This set has a dense, countable subset $L_x \subseteq \LC(\sigma)$.
 Now let $\mathcal{S}_Y$ be the set of rotational components that contain a linear approach $[\gamma]\in \LC(\sigma)$ such that there exist a representative $\gamma$ and $t>0$ with $\gamma(t)\in Y$ and $d(\gamma(t),\sigma)=t$.
 With $\mathcal{S}$ as in the proof of \autoref{prop_rotational_components_of_length_at_least_pi_are_dense}, we have
 \begin{equation*}
  \cup_{x\in Y} L_x \subseteq \mathcal{S}_Y \subseteq \mathcal{S} \subseteq \LT(\sigma) .
 \end{equation*}
 By definition, the countable set $\cup_{x\in Y} L_x$ is dense in $\mathcal{S}_Y$ and we also know from the previous proposition that the set $\mathcal{S}$ is dense in $\LT(\sigma)$. It remains to show that $\mathcal{S}_Y$ is dense in $\mathcal{S}$.

 By definition, for any $c\in\mathcal{S}$, there exists a $[\gamma]\in c$ and $t>0$ such that $d(\gamma(t),\sigma)=t$. Choose a flat local coordinate system such that $\gamma(t)$ has coordinate $(0,0)$ and $\gamma(s)=(0,s-t)$, and let $\{x_i\}$ be a sequence in $Y \cap\{(x,y):y<0,|x|<e^{1/y}\}$ converging to $\gamma(t)= (0,0)$ and $\{\gamma_i\}$ the sequence of geodesic segments from $\sigma$ to $x_i$. Then the sequence $\{[\gamma_i]\}$ of linear approaches converges to $[\gamma]$ in $\LC(\sigma)$, hence $\{\overline{[\gamma_i]}\}\subset\mathcal{S}_Y$ converges to $c$.
\end{proof}
\end{cor}

\section{What kind of topologies are possible for \texorpdfstring{$\LT(\sigma)$}{the space of rotational components}?}\label{sec:finite_top}

This section is devoted to the proof of \autoref{thm:finite_top} asserting that any finite topological space is the space of rotational components of a translation surface.

So far, we have already seen in \autoref{sec:examples} that $\LT(\sigma)$ can be different topological spaces of cardinality $1$ or $2$ (see Examples \ref{exa:harmonic_series_decoration}, \ref{exa:geometric_series_decoration} (i), and \ref{exa:geometric_series_decoration} (ii)).
However, the space of rotational components does not have to be finite. For example, the singularity of the star decoration in \autoref{exa:star_decoration} (ii) has a space of rotational components which is homeomorphic to $S^1$.
As in this example of the star decoration, the space $\LT(\sigma)$ can be Hausdorff but it does not have to be Hausdorff (cf.\ \autoref{exa:chamanara_surface} and \autoref{exa:star_decoration} (i)).

We will show now that there is indeed a large variety in spaces that can occur as a space of rotational components.

\pagebreak[3]

\begin{thm}[Every finite space occurs as $\LT(\sigma)$]\label{thm:finite_top}
 Let $Y$ be a non-empty topological space of finite cardinality. Then there exists a compact translation surface $(X, \mathcal{A})$ with a wild singularity $\sigma$ so that the space $\LT(\sigma)$ of rotational components carries the same topology as $Y$.

\begin{proof}
 Topologies on a finite set are in one-to-one correspondence with preorders “$\leq$” defined by $x\leq y$ if and only if $x$ is in the closure of $\{y\}$. So, let $n\geq 1$ and $Y = \{y_1, \ldots, y_n\}$ be a topological space of cardinality $n$ and “$\leq$” the corresponding preorder.

 We construct the translation surface $(X, \mathcal{A})$ explicitly starting with a torus with $n$ shrinking star decorations. Then we glue slits to additional tori to identify the singularities of all different shrinking star decorations to the desired singularity $\sigma$. We obtain $n$ rotational components. In a last step, we introduce more gluings for the rotational components that shall be comparable under the preorder to obtain that $\LT(\sigma)$ is homeomorphic to $Y$.

 At the beginning, let us state the setting:
 \begin{itemize}
  \item Consider a flat torus with sufficiently large injectivity radius and with $n$ copies of the shrinking star decoration as in \autoref{exa:shrinking_star_decoration}, labeled with $1,2,\dots, n$. The $i$th shrinking star decoration corresponds to $y_i \in Y$.

  \item For now, the shrinking star decoration labeled $i$ has one wild singularity $\sigma_i$ with one rotational component $c_i$ which is isometric to $[0,\infty)$. The linear approach that is corresponding to $0$ is called $[\gamma_i]$ with representative $\gamma_i$.
  
 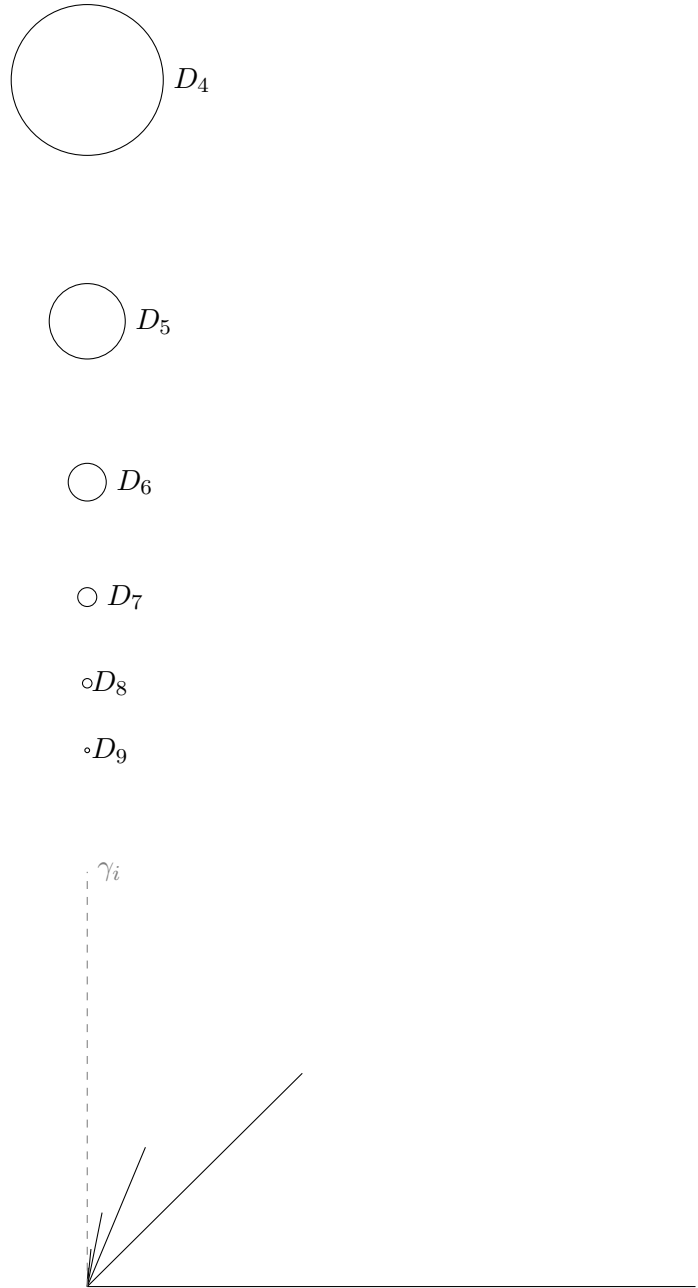
\begin{figure}
 \begin{center}
  \begin{tikzpicture}[scale=1] 
  \draw (0:0) -- (0:8);
  \draw (0:0) -- (45:4);
  \draw (0:0) -- (67.5:2);
  \draw (0:0) -- (78.75:1);
  \draw (0:0) -- (84.375:0.5);
  \draw (0:0) -- (87.1875:0.25);
  \draw[color=gray,dashed] (0,0) -- (0,5.5) node[right] {$\gamma_i$};

  \draw (0,16) circle (1);
  \path (1,16) node[right] {$D_4$};
  \draw (0,12.8) circle (0.5);
  \path (0.5,12.8) node[right] {$D_5$};
  \draw (0,32/3) circle (0.25);
  \path (0.25,32/3) node[right] {$D_6$};
  \draw (0,64/7) circle (0.125);
  \path (0.125,64/7) node[right] {$D_7$};
  \draw (0,8) circle (0.0625);
  \path (0.3,8) node {$D_8$};
  \draw (0,64/9) circle (0.03125);
  \path (0.3,64/9) node {$D_9$};
 \end{tikzpicture}
 \label{fig_shrinking_star_decoration_with_disks}
 \caption{The $i$th shrinking star decoration with additional disks. The $k$th disk is denoted $D_k$.}
 \end{center}
\end{figure}

  \item For the $i$th copy of the shrinking star decoration and every $k\geq 1$, let the $k$th disk be the disk with center $\gamma_i(\frac{4}{k})$ and radius $2^{-k}$ as shown in \autoref{fig_shrinking_star_decoration_with_disks}. This gives us an infinite sequence of disjoint disks for every shrinking star decoration.

  \item For the $k$th branch of each copy of the shrinking star decoration, we glue in a cross as in \autoref{fig:cross_with_gluings} instead of the rectangle. The arms of the cross are each glued in a similar way as the rectangle before.
  By the \emph{slit in the $k$th cross}, we mean the geodesic segment which corresponds to the top and bottom of the horizontal arms of the cross (indicated by a dashed line in \autoref{fig:cross_with_gluings}).
  We will cut open and reglue along some of these slits later.

  \item Let $T_k$ be a torus, glued from a square of side length $2^{-k}$ for each $k\geq 1$.

  \item For given $i,j \in \{1,\ldots, n\}$ define $a_k^{i,j}=k n^2+in+j$ for every $k \geq 1$. This gives $n^2$ infinite sequences in $\mathbb{N} \setminus \{0\}$ so that no two sequences have a common element.
 \end{itemize}
 
 \begin{figure}
 \begin{center}
 \begin{tikzpicture}[scale=0.85]
  \draw[-] (-4.25,0.3) -- node[above, text height=7pt]{$a$} (-2.25,0.3);
  \fill (-2.25,0.3) circle (1pt);
  \draw[-] (-2.25,0.3) -- node[above, text height=7pt]{$b$} (-0.25,0.3);
  \draw[-](-4.25,0.3) -- (-4.25,-0.3);
  \draw[-] (-4.25,-0.3) -- node[below, text height=7pt]{$b$} (-2.25,-0.3);
  \fill (-2.25,-0.3) circle (1pt);
  \draw[-] (-2.25,-0.3) -- node[below, text height=7pt]{$a$} (-0.25,-0.3);
  \draw[-](-0.25,-0.3) -- node[left, text height=7pt]{$c$} (-0.25,-2.3);
  \fill (-0.25,-2.3) circle (1pt);
  \draw[-](-0.25,-2.3) -- node[left, text height=7pt]{$d$} (-0.25,-4.3);
  \draw[dashed](-0.25,-4.3)--(0.25,-4.3);
  \draw[-](0.25,-0.3) -- node[right, text height=7pt]{$d$} (0.25,-2.3);
  \fill (0.25,-2.3) circle (1pt);
  \draw[-](0.25,-2.3) -- node[right, text height=7pt]{$c$} (0.25,-4.3);
  \draw[-] (4.25,0.3) -- node[above, text height=7pt]{$e$} (2.25,0.3);
  \fill (2.25,0.3) circle (1pt);
  \draw[-] (2.25,0.3) -- node[above, text height=7pt]{$f$} (0.25,0.3);
  \draw[-](4.25,0.3) -- (4.25,-0.3);
  \draw[-] (4.25,-0.3) -- node[below, text height=7pt]{$f$} (2.25,-0.3);
  \fill (2.25,-0.3) circle (1pt);
  \draw[-] (2.25,-0.3) -- node[below, text height=7pt]{$e$} (0.25,-0.3);
  \draw[-](-0.25,0.3) -- node[left, text height=7pt]{$h$} (-0.25,2.3);
  \fill (-0.25,2.3) circle (1pt);
  \draw[-](-0.25,2.3) -- node[left, text height=7pt]{$g$} (-0.25,4.3);
  \draw[dashed](-0.25,4.3)--(0.25,4.3);
  \draw[-](0.25,0.3) -- node[right, text height=7pt]{$g$} (0.25,2.3);
  \fill (0.25,2.3) circle (1pt);
  \draw[-](0.25,2.3) -- node[right, text height=7pt]{$h$} (0.25,4.3);
 \end{tikzpicture}
 \label{fig:cross_with_gluings}
 \caption{All segments labeled with a letter are of length $1$ whereas the other four segments are of length $2^k$.
 Segments with the same letters are glued. The left and the right side are glued to the $k$th branch of a copy of the shrinking star decoration. The top and the bottom are glued together and define the slit in the $k$th cross.}
 \end{center}
 \end{figure}
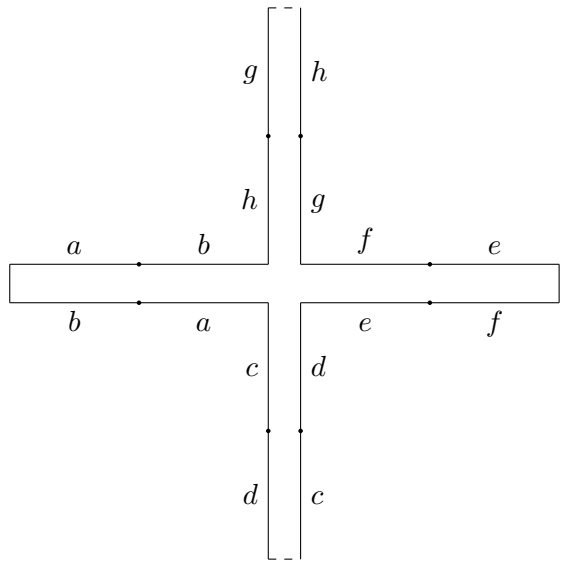

 Now we glue the copies of the shrinking star decoration in such a way that all singularities $\sigma_i$ are identified to one singularity $\sigma$. For this recall that for every copy of the shrinking star decoration and every choice of $i,j \in \{1,\ldots,n\}$, there is an infinite sequence of disks, consisting of the $a_k^{i,j}$th disk for all $k \geq 1$. Similarly, there is also the infinite sequence of the slits in the $a_k^{i,j}$th cross for all $k \geq 1$. For every singularity~$\sigma_i$, we use the sequence of the slits in the $a_k^{i,i}$th crosses of the $i$th shrinking star decoration to reduce the distance to the prospective singularity $\sigma$ to $0$. Furthermore, if the closure of the corresponding rotational component $c_i$ shall contain another rotational component~$c_j$ then we use the sequence of the slits in the $a_k^{i,j}$th crosses of the $j$th shrinking star decoration and the sequence of the $a_k^{i,j}$th disks of the $i$th shrinking star decoration. Note that this means that not all disks are used and not all slits of crosses are cut open and reglued. However, the different sequences $a_k^{i,j}$ for given $i,j$ are still useful because the disjointness of the sequences releases us from checking that no glued slits intersect.

 Let us now formally state the construction: For every $k \geq 1$ and every $i\in \{1, \ldots, n\}$, glue the slit in the $a_k^{i,i}$th cross of the $i$th shrinking star decoration to a suitable slit, i.e.\ a slit with the same holonomy vector, on the $k$th torus.
 As the side length $2^{-k}$ of the tori is shrinking, the $n$ slits for a given $k$ come the closer to each other the larger $k$ is. This makes sure that $d(\sigma_i, \sigma_j) = 0$ and hence $\sigma_i = \sigma_j$ for all $i,j \in \{1,\ldots, n\}$. The resulting singularity is called~$\sigma$. Note that the additional gluings extend the length of the rotational components but the topology of $\LT(\sigma)$ is the discrete topology for now (cf.\ the arguments in the proof of \autoref{prop:biinfinite_means_open}).

 We have to do more gluings to provide $\LT(\sigma)$ with the topology of $Y$.
 If $i\neq j$ and $y_i\leq y_j$, i.e.\ if every open set in $Y$ containing $y_i$ contains $y_j$, then we glue the slit in the $a_k^{i,j}$th cross of the $j$th shrinking star decoration to a suitable slit in the interior of the right half of the $a_k^{i,j}$th disk on the $i$th shrinking star decoration for every $k \geq 1$.
 These gluings extend the length of the rotational component $c_j$ but do not change the length of the rotational component $c_i$. More important, for every $t> 0$ and $r>0$ there exists a linear approach contained in $B(\gamma_i, t, r)$ which starts at the end point of the slit in the $a_k^{i,j}$th disk for some $k\geq 1$.
 This linear approach is contained in the rotational component~$c_j$ and hence every open neighborhood of $[\gamma_i]$ has nonempty intersection with the set of linear approaches contained in $c_j$. This means that every open neighborhood of $c_i$ in $\LT(\sigma)$ contains $c_j$.

 This ends the construction and all in all, we have shown that $\LT(\sigma)$ is homeomorphic to $Y$ as was to be proven.
 The translation surface $(X,\mathcal{A})$ that we constructed has a compact metric completion and finite area.
\end{proof}
\end{thm}

To emphasize that we can realize topological spaces of different kinds as spaces of rotational components, we use the notion of dimension (see \cite{engelking1978dimension} for definitions).

The Lebesgue covering dimension of $\LC(\sigma)$ is $1$ if $\sigma$ is a regular point or a cone angle singularity. So it is interesting that the Lebesgue covering dimension of the quotient space $\LT(\sigma)$ is not bounded if $\sigma$ is a wild singularity.

\begin{cor}[Every Lebesgue covering dimension can be realized]
 For every $n\in \mathbb{N}$, there exists a compact translation surface with one wild singularity $\sigma$ so that $\LT(\sigma)$ has Lebesgue covering dimension $n$.
 However, the small inductive dimension of the considered examples is always $1$.

\begin{proof}
 Take a set of points $p_1,\ldots,p_{n+1}, p$ with the topology characterized in the following way: a set is open if and only if it contains $p$. We can realize this as a space $\LT(\sigma)$.

 When we choose the open sets $\{p_1,p\}, \ldots, \{p_{n+1}, p\}$ as an open cover, there is no refinement of this cover and $p$ is contained in $n+1$ sets. So the Lebesgue covering dimension of $\LT(\sigma)$ is $n$.
\end{proof}
\end{cor}

Similarly to the construction in the proof of \autoref{thm:finite_top}, we can build an example where the space $\LT(\sigma)$ of rotational components is not locally compact.

\begin{exa}[A translation surface with non locally compact $\LT(\sigma)$]\label{exa:non-locally_compact}
 Consider a rooted tree of countably infinite valence. We label the vertices in the following way. The root shall be labeled by $1$. For each vertex with number $n$, the children shall be labeled by~$n\cdot q$ where $q$ runs through all primes that are greater than all prime factors of $n$.
 
 Now we build a translation surface by decorating a plane with a shrinking star for each vertex. Note that the shrinking star decorations as well as the vertices of the tree are labeled by the set $N$ of natural numbers with distinct prime factors.
 As before, we consider an infinite sequence of disks and an infinite sequence of slits in the crosses for every shrinking star decoration.
 
 Again, we glue the copies of the shrinking star decoration in such a way that all singularities are identified to one singularity $\sigma$.
 In this example, we use infinitely many additional planes instead of infinitely many additional tori.
 For every $k\geq 2$ and every $n \in N$ we glue the slit in the $2^k$th cross of the shrinking star construction with number~$n$ to a suitable slit on the $k$th plane. This is done in a way so that for fixed $n_1, n_2 \in N$ the slits in the planes have shrinking distance. This ensures that all singularities are identified to one singularity $\sigma$ with infinitely many rotational components.
 
 Furthermore, for every shrinking star decoration with number $n\in N$ and for every prime number~$q$ which is greater than all prime factors of $n$, we glue a suitable slit in the right half of the $q$th disk to the slit in the $q$th cross of the shrinking star decoration with number $nq$. Note that the star decoration with number $nq$ corresponds to a child of the vertex that corresponds to the star decoration with number $n$. 
 
 To study the topology of the space of rotational components, we identify each vertex of the tree with its corresponding shrinking star decoration and its corresponding rotational component for convenience.
 As a consequence of the gluings, every neighborhood of a vertex contains all but finitely many of its children. In particular, any neighborhood of the root must contain a tree of infinite height. Let $U$ be such a neighborhood of the root and $\{v_i\}_{i \in \mathbb{N}} \subset U$ an infinite sequence of vertices where $v_{i+1}$ is a child of $v_i$. Note that the subspaces $U_k=U \setminus \{v_i : i>k\}$ are open in the subspace topology on $U$. Then $\{U_k\}$ is an open cover of $U$ without finite subcover, hence $U$ cannot be compact. Therefore, the root and its corresponding rotational component do not have a compact neighborhood.
\end{exa}

\section{Topology of \texorpdfstring{$\LT(\sigma)$}{the space of rotational components} does not determine topology of \texorpdfstring{$\LC(\sigma)$}{the space of linear approaches}}

When we pass from the space $\LC(\sigma)$ of linear approaches to the quotient space $\LT(\sigma)$ of rotational components, we apparently lose information. It is a natural question to raise how much information we lose.
In \autoref{thm:info_loss} we show that there are uncountably many translation surfaces with non-homeomorphic spaces of linear approaches but with homeomorphic spaces of rotational components.

\pagebreak

\begin{thm}[$\LT(\sigma)$ does not determine $\LC(\sigma)$]\label{thm:info_loss}
 There are uncountably many translation surfaces $(X_r, \mathcal{A}_r)$ with exactly one singularity~$\sigma_r$, with pairwise non-homeomorphic spaces of linear approaches $\LC(\sigma_r)$, but with homeomorphic spaces of rotational components $\LT(\sigma_r)$.

\begin{proof}
 Let $\mathbb{K} = \{0, 1\}^\mathbb{N}$ and consider the equivalence relation on $\mathbb{K}$ defined by $r \sim r'$ if and only if some shifts of $r$ and $r'$ are equal, i.e.\ if the tails of $r$ and $r'$ agree up to shift.

 We define a translation surface $(X_r, \mathcal{A}_r)$ for every $r\in \mathbb{K}$ but not in the equivalence class of $(0,0,0,\ldots)$ by modifying the stack of boxes from \autoref{exa:stack_of_boxes}. The translation surface $(X_r, \mathcal{A}_r)$ will have exactly one singularity $\sigma_r$ and this singularity will have one rotational component.
 Then we show that if $\LC(\sigma_r)$ and $\LC(\sigma_{r'})$ are homeomorphic then $r$ and $r'$ are in the same equivalence class.
 This proves the statement since there are uncountably many equivalence classes in $\mathbb{K}$.
 
 Let $r=(r_n)\in \mathbb{K}$ but not in the equivalence class of $(0,0,0,\ldots)$.
 Recall the stack of boxes from \autoref{exa:stack_of_boxes} and choose the two sequences $H$ and $W$ by $h_n = 1$ and $w_n = 2^{-n}$ for every $n\in \mathbb{N}$.
 We modify the stack of boxes by gluing in additional rectangles, depending on $r$:
 Consider vertical geodesic segments $v_n$ of length $2^{-2n}$ starting at each $A_n$ for $n\geq 2$ (see \autoref{fig_stack_of_boxes_with_vertical_line_segments}).
 For each $n \geq 2$ with $r_n=1$, we cut open along $v_n$ and glue the two sides of the slit $v_n$ to the two vertical sides of a rectangle of width $1$, then we identify the top and bottom of the rectangle.
 We call the resulting translation surface $(X_r, \mathcal{A}_r)$.
 As in \autoref{exa:stack_of_boxes}, it has one singularity $\sigma_r$ and one rotational component isometric to $(0,\infty)$.
 Hence the spaces $\LT(\sigma_r)$ are homeomorphic for all $r$.

 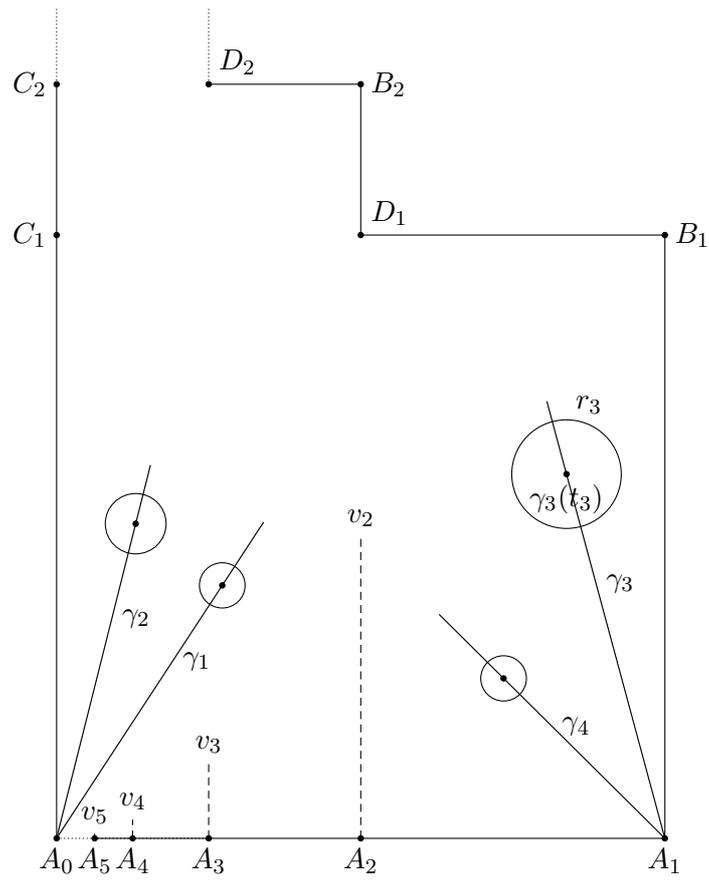
\begin{figure}
 \begin{center}
 \begin{tikzpicture}[scale=1]
  \draw[densely dotted] (0,0) -- (2,0); 
  \draw (.5,0) -- (1,0) -- (2,0) -- (4,0) -- (8,0) -- (8,8) -- (4,8) -- (4,10) -- (2,10);
  \draw (0,0) -- (0,10);
  \draw[densely dotted] (0,10) -- (0,11);
  \draw[densely dotted] (2,10) -- (2,11);

  \draw[fill] (0,0) circle (1pt) node[below] {$A_0$}; 
  \draw[fill] (0.5,0) circle (1pt) node[below] {$A_5$};
  \draw[densely dashed](0.5,0)--(0.5,0.0625) node[pos=1,above]{$v_5$};
  \draw[fill] (1,0) circle (1pt) node[below] {$A_4$};
  \draw[densely dashed](1,0)--(1,0.25) node[pos=1,above]{$v_4$};
  \draw[fill] (2,0) circle (1pt) node[below] {$A_3$};
  \draw[densely dashed](2,0)--(2,1) node[pos=1,above]{$v_3$};
  \draw[fill] (4,0) circle (1pt) node[below] {$A_2$};
  \draw[densely dashed](4,0)--(4,4) node[pos=1,above]{$v_2$};
  \draw[fill] (8,0) circle (1pt) node[below] {$A_1$};
  \draw[fill] (8,8) circle (1pt) node[right] {$B_1$};
  \draw[fill] (4,10) circle (1pt) node[right] {$B_2$};
  \draw[fill] (0,8) circle (1pt) node[left] {$C_1$};
  \draw[fill] (0,10) circle (1pt) node[left] {$C_2$};
  \draw[fill] (4,8) circle (1pt) node[above right] {$D_1$};
  \draw[fill] (2,10) circle (1pt) node[above right] {$D_2$};

  \draw[fill] (0,0) -- node[right, pos=0.7]{$\gamma_1$} ++(57:4) circle (1pt) -- ++(57:1);
  \draw (57:4) circle (0.3);
  
  \draw[fill] (0,0) -- node[right, pos=0.7]{$\gamma_2$} ++(76:4.3) circle (1pt) -- ++(76:0.8);
  \draw (76:4.3) circle (0.4);
  
  \draw[fill] (8,0) -- node[right, pos=0.7]{$\gamma_3$} ++(105:5) circle (1pt) node[below] {$\gamma_3(t_3)$} -- ++(105:1);
  \draw (8,0) ++ (105:5) circle (0.72);
  \draw (7,5.75) node {$r_3$};
  
  \draw[fill] (8,0) -- node[right, pos=0.7]{$\gamma_4$} ++(135:3) circle (1pt) -- ++(135:1.2);
  \draw (8,0) ++ (135:3) circle (0.3);
 \end{tikzpicture}
 \label{fig_stack_of_boxes_with_vertical_line_segments}
 \caption{Stack of boxes with additional vertical slits (not to scale in vertical direction).}
 \end{center}
\end{figure}

 However, each geodesic segment starting at $A_0$ or $A_1$ defines a linear approach whose neighborhoods change drastically with each slit $v_n$ that is glued to the sides of a rectangle. We use this to elaborate tools to recover the equivalence class of $r$ intrinsically from the topology of $\LC(\sigma_r)$ for an $r\in \mathbb{K}$.
 Let $c$ be the unique rotational component of $\sigma_r$ with the one-dimensional translation structure as in \autoref{def:translation_struct}. Then there exists an isometry $(0, \infty) \to c$. As $c$ is the only rotational component of $\sigma_r$, we also have a bijective map $(0,\infty) \to \LC(\sigma_r)$. Let $\mathbb{R}_+^r$ be the set $(0,\infty)$, equipped with the pullback of the topology of $\LC(\sigma_r)$. 
 In the following, we use the parameter $x \in \mathbb{R}_+^r$ to represent linear approaches in $\LC(\sigma_r)$ and we let $x$ run through $\mathbb{R}_+^r$, i.e.\ the linear approach is rotating counterclockwise around the singularity.

 To record when the rotated linear approach is contained in an additional cylinder, we use four distinguished neighborhoods in $\LC(\sigma_r)$ (see \autoref{fig_stack_of_boxes_with_vertical_line_segments}).
 Let
 \begin{equation*}
  \mathcal{U}_r = B( \gamma_1, t_1, \rho_1), \
  \mathcal{U'}_r = B( \gamma_2, t_2, \rho_2), \
  \mathcal{V}_r = B( \gamma_3, t_3, \rho_3), \
  \mathcal{V'}_r = B( \gamma_4, t_4, \rho_4), \
 \end{equation*}
 where for $i \in \{1,\ldots, 4\}$, $\gamma_i$ is a representative of a linear approach that is represented by some $x_i$ with $0<x_1<x_2< \frac{\pi}{2} < x_3<x_4<\pi$ and $t_i > 0$, $\rho_i > 0$ are small enough so that the $B(\gamma_i (t_i), \rho_i)$ are disjoint from all slits $v_n$, from the edges of the boxes, and from each other.

 Now let $n\geq 2$ be large enough and let $y_n$ be the number which corresponds to the horizontal linear approach starting at $A_n$ and going to the right.
 For $y\in \mathbb{R}_+^r$ with $y>y_n$ but sufficiently close to $y_n$, $y$ is in neither of $\mathcal{U}_r$, $\mathcal{U'}_r$, $\mathcal{V}_r$, or $\mathcal{V'}_r$.
 When $y$ is increased and gets closer to $y_n + x_1$ we eventually have $y \in \mathcal{U}_r$. When $y$ increases further, it leaves $\mathcal{U}_r$ to enter $\mathcal{U'}_r$.
 Now there are two cases to consider.

 Case 1: If $r_n = 1$, when increasing $y$ further the linear approach will eventually be contained in the cylinder attached to $v_n$ (starting at the lower right corner of the glued-in rectangle),
 and after additional time $\pi$ will not be contained in the cylinder any more but starting at the tip of the slit $v_n$.
 After running for about another $\frac{\pi}{2} +x_1$, $y$ will reenter $\mathcal{U}_r$.
 Then it will enter $\mathcal{U'}_r$, $\mathcal{V}_r$, $\mathcal{V'}_r$, to finally enter the same cylinder attached to $v_n$ from the left side.
 After running for another $\pi$, $y$ will exit the cylinder and enter $\mathcal{V}_r$ and $\mathcal{V'}_r$ again.
 When $y$ reaches $y_n+5\pi$, the linear approach will leave the first box and enter the $n$th box from the upper right corner at $B_n$.
 
 Case 2: If $r_n = 0$, after going through $\mathcal{U}_r$ and $\mathcal{U'}_r$, $y$ will simply go through $\mathcal{V}_r$ and $\mathcal{V'}_r$ and then starts at $B_n$ when reaching $y_n+\pi$.

 By recording the successive passages in each of the four neighborhoods, writing $1$ when seeing $(\mathcal{U}_r, \mathcal{U'}_r, \mathcal{U}_r, \mathcal{U'}_r, \mathcal{V}_r, \mathcal{V'}_r, \mathcal{V}_r, \mathcal{V'}_r)$ and $0$ when seeing $(\mathcal{U}_r, \mathcal{U'}_r, \mathcal{V}_r, \mathcal{V'}_r)$, we record a sequence in $\mathbb{K}$ that eventually agrees with $r$ up to shift. Note that the exact values of $0<x_1<x_2<\frac{\pi}{2}<x_3<x_4<\pi$ (and also of the $t_1,\ldots,t_4$ and $\rho_1,\ldots,\rho_4$) do not matter for this argument. We could even choose neighborhoods of $[\gamma_1], \ldots, [\gamma_4]$ that are not of the same form as the elements of the subbasis defined in \autoref{lem:subbasis_LCpoint}. In the latter case, we have to be careful to record the neighborhoods only up to multiple leavings and reentries in one neighborhood before we enter the next neighborhood.

 To finish the proof, let $\varphi$ be a homeomorphism between $\LC(\sigma_r)$ and $\LC(\sigma_{r'})$. Then $\varphi$ induces a homeomorphism $f \colon \mathbb{R}_+^r \rightarrow \mathbb{R}_+^{r'}$.
 We want to show that the tails of $r$ and $r'$ agree up to shift by looking at the recorded sequences. Hence we have to examine the images of the sets $\mathcal{U}_r, \mathcal{U'}_r, \mathcal{V}_r, \mathcal{V'}_r$ under $\varphi$.

 By \autoref{lem:neighborhood_linear_approaches} we have that all elements of $\LC(\sigma_r)$ or $\mathbb{R}_+^r$ in $(\pi,\infty)$ have a path-connected neighborhood. On the other hand, a close look reveals that the elements in $(0,\pi]$ do not have arbitrarily small path-connected neighborhoods. As a homeomorphism must preserve the set of points that have path-connected neighborhoods, we have $f((0,\pi]) = (0,\pi]$.
 Furthermore, we can distinguish the neighborhoods of $\frac{\pi}{2}$ in $\mathbb{R}_+^r$ and $\mathbb{R}_+^{r'}$ from the neighborhoods of other points in $(0, \pi]$: The linear approach corresponding to $\frac{\pi}{2}$ has the same direction as the slits. Then again the lengths of the linear approaches defined by the slits tend to $0$. This implies that the linear approaches corresponding to the slits cannot be contained in a small neighborhood of the linear approach corresponding to $\frac{\pi}{2}$. This makes the neighborhoods of $\frac{\pi}{2}$ in $\mathbb{R}_+^r$ and $\mathbb{R}_+^{r'}$ distinguishable from the neighborhoods of other points.
 Hence, $f(\frac{\pi}{2})= \frac{\pi}{2}$.

 We show next that the map $f \colon \mathbb{R}_+^r \to \mathbb{R}_+^{r'}$ is monotone.
 Let $I_r$ be an open interval in $(0,\infty)$ of length shorter than $\pi$ which does not end at $0$ and so that the open interval $f(I_r)$ is also of length shorter than $\pi$. Equip $I_r$ with the subspace topology from $\mathbb{R}_+^r$.
 We show that this subspace topology is the same as the Euclidean topology of~$I_r$.
 For this, choose a linear approach $[\gamma]$ representing a point in $I_r$.
 It does not lie at an endpoint of $\mathbb{R}_+^r$ resp.\ of the rotational component, equipped with the translation structure with boundary as in \autoref{def:translation_struct}.
 Hence, given a representative $\gamma$ of $[\gamma]$ and a $t>0$ such that $\gamma(t)$ is defined, $\gamma((0,t])$ lies in the image of an angular sector centered at $\sigma_r$.
 This means that any neighborhood of $[\gamma]$ in $I_r$ under the subspace topology must contain a neighborhood under the Euclidean topology.
 On the other hand, there is an embedded open sector centered at $\sigma_r$ such that all linear approaches on it are represented by points on $I_r$, and any open subinterval of $I_r$ corresponds to an open sub-sector, hence must be open in the subspace topology.
 The map $f_{|I_r} \colon I_r \to \mathbb{R}_+^{r'}$ is homeomorphic to its image and so it is monotone due to the argument above. This means that $f \colon \mathbb{R}_+^{r} \to \mathbb{R}_+^{r'}$ is locally monotone and so it is monotone.

 As a consequence of the monotonicity, we can choose $x_1,x_2,x_3,x_4$ with $0<x_1<x_2<\frac{\pi}{2}<x_3<x_4<\pi$ and $0<f(x_1)<f(x_2)<\frac{\pi}{2}<f(x_3)<f(x_4)<\pi$.

 When choosing the $\rho_i$ sufficiently small, the discussion above shows that the recorded sequences for $y$ in $\LC(\sigma_r)$ and in $\LC(\sigma_{r'})$ consisting of $0$'s and $1$'s agree up to shift, hence $r$ and $r'$ agree up to shift and so they are in the same equivalence class, as was to be shown.
\end{proof}
\end{thm}

\bibliographystyle{alpha} 
\bibliography{sing}

\end{document}